\documentclass[aos,MSNbibl,seceqn,nameyear,dvips]{arximspdf}

\doi{10.1214/12-AOS975} 
\volume{40}
\issue{2}
\pubyear{2012}
\firstpage{727}
\lastpage{758}

\makeatletter

\newcommand{\cal}{\mathcal}

\newtheorem{theorem}{Theorem}[section]
\newtheorem{cor}{Corrolary}[section]
\newtheorem{prop}{Proposition}[section]

\newproclaim{Assumption}{Assumption}

\newcommand{\ep}{\epsilon}
\newcommand{\Lam}{\Lambda}
\newcommand{\al}{\alpha}
\newcommand{\si}{\sigma}

\makeatother

\begin{document}
\begin{frontmatter}

\title{A specification test for nonlinear nonstationary~models}
\runtitle{Specification test for nonstationary model}

\begin{aug}
\author[A]{\fnms{Qiying} \snm{Wang}\corref{}\thanksref{t1}\ead[label=e1]{qiying@maths.usyd.edu.au}}
\and
\author[B]{\fnms{Peter C. B.} \snm{Phillips}\thanksref{t2}\ead[label=e2]{peter.phillips@yale.edu}}
\runauthor{Q. Wang and P. C. B. Phillips}
\affiliation{University of Sydney,
and Yale University, University of Auckland, University~of~Southampton,
Singapore Management University}
\address[A]{School of Mathematics and Statistics\\
University of Sydney\\
NSW 2006\\
Australia\\
\printead{e1}}
\address[B]{Department of Economics\\
Yale University\\
30 Hillhouse Avenue\\
New Haven, Connecticut 06520\\
USA\\
\printead{e2}} 
\end{aug}

\thankstext{t1}{Supported by the Australian Research Council.}

\thankstext{t2}{Supported by the NSF Grant SES 09-56687.}

\received{\smonth{10} \syear{2010}}
\revised{\smonth{1} \syear{2012}}

%
\begin{abstract}
We provide a limit theory for a general class of kernel smoothed
U-statistics that may be used for specification testing in time series
regression with nonstationary data. The test framework allows for
linear and nonlinear models with endogenous regressors that have
autoregressive unit roots or near unit roots. The limit theory for the
specification test depends on the self-intersection local time of a
Gaussian process. A new weak convergence result is developed for
certain partial sums of functions involving nonstationary time series
that converges to the intersection local time process. This result is
of independent interest and is useful in other applications.
Simulations examine the finite sample performance of the test.
\end{abstract}

%

\end{frontmatter}

\section{Introduction}\label{sec1}

One of the advantages of nonparametric modeling is the opportunity for
specification testing of particular parametric models against general
alternatives. The past three decades have witnessed many developments
in such specification tests involving nonparametric and semiparametric
techniques that allow for independent, short memory and long-range
dependent data. Recent research on the nonparametric modeling of
nonstationary data opens up some new possibilities that seem relevant
to applications in many fields, including nonlinear diffusion models in
continuous time [Bandi and Phillips (\citeyear{BanPhi03},
\citeyear{BanPhi07})] and cointegration models in economics and
finance.

Cointegration models were originally developed in a linear parametric
framework 
that has been widely
used in econometric applications. That framework was extended in Park
and Phillips (\citeyear{ParPhi99}, \citeyear{ParPhi01}) to allow for nonlinear parametric formulations under certain
restrictions on the function nonlinearity. While considerably
broadening the
class of allowable nonstationary models, the potential for parametric
misspecification in these models is still present and is important to test
in applied work.

The hypothesis of linear cointegration is of particular interest in
this context, given the vast empirical literature. Recent papers by
\citet{KarMykTjs07}, Wang and Phillips
(\citeyear{WanPhi09N1}, \citeyear{WanPhi09N2}, \citeyear{WanPhi11}) and
\citet{Sch} have developed asymptotic theory for nonparametric
kernel regression of nonlinear nonstationary systems. This work
facilitates the comparison of various parametric specifications against
a more general nonparametric nonlinear alternative. Such comparisons
may be based on weighted sums of squared differences between the
parametric and nonparametric estimates of the system or on a~kernel-based
U-statistic test which uses\break a~smoothed version of the
parametric estimator in its construction [e.g., \citet{Gao07},
Chapter 3].

A major obstacle in the development of such specification tests is the
technical difficulty of developing a limit theory for these weighted
sums which typically involve kernel functions with multiple
nonstationary regressor arguments. Few results are currently available,
and because of this shortage, attempts to develop specification tests
for nonlinear regression models with nonstationarity have been highly
specific and do not involve nonparametric alternatives or kernel
methods. Some examples of recent work in parametric models include Choi
and Saikonnen (\citeyear{ChoSai04}, \citeyear{ChoSai10}),
\citet{Mar08}, \citet{HonPhi10} and \citet{KasPhi}. An
exception is the recent work for testing linearity in autoregression
and parametric time series regression by Gao et~al.
(\citeyear{Gaoetal09N1}, \citeyear{Gaoetal09N2}) who obtained a limit distribution
theory for a kernel based specification test in a setting that involves
martingale difference errors and random walk regressors.

The present paper makes a related contribution and seeks to provide
a~general theory of specification tests that is applicable for a
wider class of nonstationary regressors that includes both unit root
and near unit root processes. The latter are important in practical
work where a unit root restriction is deemed too restrictive. The
paper contributes to this emerging literature in two ways. First, we
provide a limit theory for a general class of kernel-based
specification tests of parametric nonlinear regression models that
allows for
near unit root processes driven by short memory (linear process)
errors. This limit theory should be widely applicable to
specification testing in nonlinear cointegrated systems.

Second, the limit theory of the specification test involves the
self-intersec\-tion local time of a Gaussian limit process. The result
requires establishing weak convergence to this self-intersec\-tion local time
process, which is of independent interest, and a feasible central limit
theorem involving an empirical estimator of the intersection local time that
can be used to construct the test statistic. Thus, the results provide some
new theories for intersection local time, weak convergence and specification
test asymptotics that are relevant in applications.

The paper is organized as follows. Section \ref{sec2} lays out the
nonparametric and parametric models and assumptions. Section \ref{sec3} gives
the main results on specification test limit theory. Section \ref{sec4}
reports some simulation evidence on test performance. Section \ref{sec5}
provides the weak convergence theory for intersection local time.
Section \ref{sec6} gives proofs of the main theorems in Section \ref{sec3}. The proofs
of the local time limit theory in Section \ref{sec5} and some supplemental
technical results in Section~\ref{sec6} can be found
in the supplementary material [\citet{WanPhi}].

\section{Model and assumptions}\label{sec2}

We consider the nonlinear cointegrating regression model
%
%
\begin{equation} \label{mo1}
y_{t+1}=f(x_{t})+u_{t+1},\qquad t=1,2,\ldots,n,
\end{equation}
where $u_{t}$ is a stationary error process, and $x_{t}$ is a nonstationary
regressor. We are interested in testing the null hypothesis
\[
H_{0}\dvtx f(x)=f(x,\theta), \qquad\theta\in\Omega_{0},
\]
for $x\in R$, where $f(x,\theta)$ is a given real function indexed by a
vector $\theta$ of unknown parameters which lie in the parameter space
$%
\Omega_{0}$.

To test $H_{0}$ we make use of the following kernel-smoothed test statistic:
%
%
\begin{equation} \label{main1}
S_{n}=\sum_{s,t=1,s\not=t}^{n}\hat{u}_{t+1}\hat{u}_{s+1} K[%
(x_{t}-x_{s})/h],
\end{equation}
involving the parametric regression residuals $\hat
{u}_{t+1}=y_{t+1}-f(x_{t},%
\hat{\theta})$, where $K(x)$ is a nonnegative real kernel function,
$h$ is
a bandwidth satisfying $h\equiv h_{n}\rightarrow0$ as the sample size $
n\rightarrow\infty$ and $\hat{\theta}$ is a parametric estimator of
$\theta
$ under the null $H_0$, that is consistent whenever $\theta\in\Omega_0$.

The statistic $S_{n}$ in (\ref{main1}) has commonly been applied to
test parametric specifications in stationary time series regression
[see \citet{Gao07}] and was used by Gao et~al.
(\citeyear{Gaoetal09N1}, \citeyear{Gaoetal09N2}) to test for linearity
in autoregression and a parametric conditional mean function in time
series regression involving a~random walk regressor. $S_{n}$ is a
weighted U-statistic with kernel weights that depend on standardized
differentials $(x_{t}-x_{s})/h$ of the regressor. The weights focus
attention in the statistic on those components in the sum where the
nonstationary regressor $x_{t}$ nearly intersects itself. This
smoothing scheme gives prominence to product components
$\hat{u}_{t+1}\hat{u}_{s+1}$ in the sum where~$s$ and~$t$ may differ
considerably but for which the corresponding regressor process takes
similar values (i.e., $x_{t},x_{s}\simeq x$ for some $x$), thereby
enabling a test of~$H_{0}$.\vadjust{\goodbreak}

The difficulty in the development of an asymptotic theory for
$S_{n}$ stems from the presence of the kernel weights
$K((x_{t}-x_{s})/h)$. The behavior of these weights depends on the
self intersection properties of $x_{t}$ in the sample, and, as
$n\rightarrow\infty$, this translates into the corresponding
properties of the stochastic process to which a standardized version of
$%
x_{t}$ converges. To establish asymptotics for $S_{n}$, we need to
account for this limit behavior, which leads to a new limit theory
involving the self-intersection local time of a~Gaussian process
(i.e., the local time for which the process intersects
itself).\looseness=1

We use the following assumptions in our development.
%
%
\begin{Assumption}\label{Assumption1}
(i) $\{\epsilon_{t}\}_{t\in
\mathbf{Z}}$ is a sequence of independent and identically
distributed (i.i.d.) continuous random variables with
$E\epsilon
_{0}=0$, $E\epsilon_{0}^{2}=1$, and with the
characteristic function $\varphi(t)$ of $\epsilon
_{0}$ satisfying $|t||\varphi(t)|\rightarrow0$, as
$|t|\rightarrow
\infty$.

(ii)
%
%
\begin{equation} \label{m100}\qquad
x_{t}=\rho x_{t-1}+\eta_{t},\qquad x_{0}=0,\qquad \rho=1+\kappa
/n,\qquad
1\leq t\leq n,
\end{equation}
where $\kappa$ is a constant and $\eta
_{t}=\sum_{k=0}^{\infty}\phi_{k}\epsilon_{t-k}$ with $\phi
\equiv\sum_{k=0}^{\infty}\phi_{k}\not=0$ and
$\sum_{k=0}^{\infty}k^{1+\delta}|\phi_{k}|<\infty$ for some $\delta>0$.
\end{Assumption}
%
%
\begin{Assumption}\label{Assumption2}
(i) $\{u_{t},\mathcal{F}%
_{t}\}_{t\geq1}$, where $\mathcal{F}_{t}$ is a sequence of
increasing $\sigma$-fields which is independent of $\epsilon
_{k},k\geq
t+1$, forms a martingale difference satisfying
$E(u_{t+1}^{2}\mid%
\mathcal{F}_{t})\rightarrow_{\mathrm{a.s.}}\sigma^{2}>0$ as
$t\rightarrow
\infty$ and $\sup_{t\geq1}E(|u_{t+1}|^{4}\mid\mathcal{F}
_{t})<\infty$.

(ii) $x_{t}$ is adapted to $ \mathcal
F_{t}$, and
there exists a correlated vector Brownian motion $(W,V)$
such that
%
%
\begin{equation}\label{tt1}
\Biggl(\frac{1}{\sqrt{n}}\sum_{j=1}^{[nt]}\epsilon_{j}, \frac
{1}{\sqrt{n}%
\sigma}\sum_{j=1}^{[nt]}u_{j+1}\Biggr)\Rightarrow_{D}(W(t),V(t))
\end{equation}
on $D[0,1]^{2}$ as $n\rightarrow\infty$.
\end{Assumption}
%
%
\begin{Assumption}\label{Assumption3}
$K(x)$ is a nonnegative real function satisfying\break $\sup_{x}K(x)<\infty$
and $\int K(x) \,dx<\infty$.
\end{Assumption}
%
%
\begin{Assumption}\label{Assumption4}
(i) There is a
sequence of positive real numbers
$\delta_{n}$ satisfying $\delta_{n}\rightarrow0$
as $n\rightarrow\infty$ such that ${\sup_{\theta\in\Omega
_{0}}}\|\hat{\theta}-\theta\|=o_{P}(\delta_{n})$, where
\mbox{$\|\cdot\|$} denotes the Euclidean norm.

\vspace*{-8pt}

\begin{longlist}[(iii)]
\item[(ii)] There exists
some $\varepsilon_{0}>0$ such that $\frac{\partial
^{2}f(x,t )}{\partial t ^{2}}$ is continuous in both
$x\in R$ and $t \in\Theta_0$, where $\Theta_0=\{
t \dvtx\|t -\theta\|\leq\varepsilon_{0},\theta\in
\Omega_{0}\}$.

\item[(iii)] Uniformly for $\theta\in\Omega_{0}$,
\[
\biggl| \frac{\partial f(x,t )}{\partial t }%
\bigg|_{t =\theta} \biggr| +\biggl| \frac{%
\partial^{2}f(x,t )}{\partial t ^{2}}\bigg|_{t =\theta} \biggr| \leq
C(1+|x|^{\beta})
\]
for some constants $\beta\geq0$ and $C>0$.

\item[(iv)]
Uniformly for $\theta\in\Omega_{0}$, there exist $0<\gamma'
\leq1$ and $\max\{0,3/4-2\beta\}<\gamma\le1$ such that
%
%
\begin{equation}\label{tz2}\qquad
| g(x+y,\theta)-g(x,\theta)| \leq C|y|^{\gamma}
\cases{1+|x|^{\beta-1}+|y|^{\beta}, &\quad if $\beta>0$, \cr
1+|x|^{\gamma' -1}, &\quad if $\beta=0$,}
\end{equation}
for any $x,y\in R$, where $g(x,t )=\frac{\partial
f(x,t )}{\partial t}$.
\end{longlist}
\end{Assumption}
%
%
\begin{Assumption}\label{Assumption5}
$nh^{2}\rightarrow\infty
$, $\delta_{n}^{2}n^{1+\beta}\sqrt{h}\rightarrow0$ and
$nh^{4}\log^{2}n\rightarrow0$, where $\beta$ and
$\delta_{n}^{2}$ are defined as in Assumption \ref{Assumption4}. Also,
$\int(1+|x|^{2\beta+1})K(x)\,dx<\infty$ and $E|\epsilon
_{0}|^{4\beta+2}<\infty$.
\end{Assumption}

Assumption \ref{Assumption1} allows for both a unit root ($\kappa=0$)
and a near unit root ($\kappa\neq0$) regressor by virtue of the
localizing coefficient $\kappa$ and is standard in the near integrated
regression framework [Phillips (\citeyear{Phi87}, \citeyear{Phi88}),
\citet{ChaWei87}]. Compared to the estimation theory developed in
Wang and Phillips (\citeyear{WanPhi09N1}, \citeyear{WanPhi09N2}) and
for technical convenience in the present work, we impose the stronger
summability condition $\sum_{k=0}^{\infty
}k^{1+\delta}|\phi_{k}|<\infty$ for some $\delta>0$ on the coefficients
of the linear process $\eta_{t}=\sum_{k=0}^{\infty
}\phi_{k}\epsilon_{t-k}$ driving the regressor $x_{t}$. Under these
conditions, it is well known that the standardized process $x_{[
nt%
] ,n}=x_{[ nt] }/\sqrt{n}\phi$ converges weakly
to the
Gaussian process $G( t) =\int_{0}^{t}e^{\kappa
(t-s)}\,dW(s)$, where
$W(t)$ is a standard Brownian motion. See (\ref{no1}) below
or \citet{PhiSol92}.

Assumption \ref{Assumption2}(i) is a standard martingale difference
condition on the
equation innovations $u_{t}$, so that
$\operatorname{cov}(u_{t+1},x_{t})=E[x_{t}E(u_{t+1}%
\mid\mathcal{F}_{t})]=0$. \citet{WanPhi09N2} allowed for
endogeneity in their nonparametric structure, so the equation error
could be serially dependent and cross-correlated with~$x_{t}$ for
$|t-s|\leq m_{0}$ for some finite $m_{0}$. It is not clear at the
moment if the results of the present paper on testing extend to the
more general error structure considered in \citet{WanPhi09N2},
but simulation results suggest that this may be so. Assumption
\ref{Assumption2}(ii) is a standard functional law for partial sum
processes [e.g., \citet{ParPhi01}].

Assumption \ref{Assumption3} is a standard condition on $K(x)$ as in
the stationary
situation. The integrability condition is weaker than the common alternative
requirement that $K(x)$ has compact support.

As seen in Assumption \ref{Assumption5}, the sequence $\delta_{n}$ in
Assumption \ref{Assumption4}(i) may be chosen as
$\delta_{n}^{2}=n^{-(1+\beta)/2}h^{-1/8}$. As $h\rightarrow0$ and
$\kappa=0$ in (2.3), Assumption \ref {Assumption4}(i) holds under very
general conditions, such as those of Theorem 5.2 in
\citet{ParPhi01}. Indeed, by \citet{ParPhi01}, we may choose
$\hat{\theta}$ such that ${\sup_{\theta\in\Omega
_{0}}}\|%
\hat{\theta}-\theta\|=O_{P}(n^{-(1+\beta)/2})$, under our
Assumption \ref{Assumption4}(ii)--(iv). Assumption \ref
{Assumption4}(ii)--(iv) is quite weak and
includes a
wide class of functions. Typical examples include polynomial forms like
$%
f(x,\theta)=\theta_{1}+\theta_{2}x+\cdots+\theta_{k}x^{k-1}$, where
$\theta
=(\theta_{1},\ldots,\theta_{k})$, power functions like $f(x,a ,b
,c )=a+b x^{c }$, shift functions like $f(x,\theta
)=x(1+\theta x)I(x\geq0)$ and weighted exponentials such as $f(x,a
,b )=(a +b e^{x})/(1+e^{x})$. However, Assumption \ref{Assumption4} excludes
models where $f(x,\theta)$ is integrable, because parametric rates of
convergence are known to be $O( n^{1/4}) $ in this case
[see \citet{ParPhi01}]. It seems that cases with integrable $f(x,\theta)$
require different techniques and these are left for future investigation.

As in estimation limit theory, the condition in Assumption \ref
{Assumption5} that
the bandwidth $h$ satisfies $nh^{2}\rightarrow\infty$ is
necessary. The further condition that $nh^{4}\log^{2}n\rightarrow
0$ restricts the choice of $h$ and, at least with the techniques
used here, seems difficult to relax in the general case studied in
the present work, although it may be substantially relaxed in less
general models as discussed later in the paper. The condition that
$\delta_{n}^{2}n^{1+\beta}\sqrt{h}\rightarrow0$
holds automatically if $\sup_{\theta\in\Omega_{0}}\|%
\hat{\theta}-\theta\|=O_{P}(n^{-(1+\beta)/2})$.
As explained above, the latter condition holds true under very general
settings such as
Assumption \ref{Assumption4}(ii)--(iv). We also impose a
higher moment condition on the innovation $\ep_0$ in Assumption~\ref
{Assumption5}
which helps in the development of
the limit theory.

\section{Main results on specification}\label{sec3}

The limit distribution of $S_{n}$ under standardization involves
nuisance parameters $\sigma$ and $\phi$, which are the limit of
$Eu_{t}^{2}$ as $t\rightarrow\infty$ and the sum of coefficients
of the linear process appearing in Assumption \ref{Assumption1}; see
Corollary \ref
{cor1} below. While convenient, this formulation obviously restricts
direct use of the result in applications. The dependence on the
nuisance parameters can be simply removed by self-normalization.
Indeed, by defining
\[
V_{n}^{2}=\sum_{s,t=1,s\not=t}^{n}\hat{u}_{t+1}^{2}\hat
{u}_{s+1}^{2} K^{2}%
[(x_{t}-x_{s})/h],
\]
we have the following main result.
%
%
\begin{theorem}
\label{th3} Under Assumptions \ref{Assumption1}--\ref{Assumption5} and
the null hypothesis, we have
%
%
\begin{equation} \label{re3}
\frac{S_{n}}{\sqrt{2} V_{n}}\rightarrow_{D}N,
\end{equation}
where $N$ is a standard normal variate.
\end{theorem}

The limit in Theorem \ref{th3} is normal and does not depend on any
nuisance parameters. As a test statistic,
$Z_{n}=S_{n}/\sqrt{2} V_{n}$ has a big advantage in applications.
In order to investigate the asymptotic power of the test, we
consider the local alternative models
\[
H_1\dvtx f(x)=f(x, \theta)+\rho_n m(x),
\]
where $\theta\in\Omega_0$, $\rho_n$ is a sequence of constants,
and $m(x)$ is a real function. This kind of local alternative model
is commonly used in the theory of nonparametric inference involving
stationary data; see, for instance, \citet{HorSpo01}.\vadjust{\goodbreak}
%
%
\begin{Assumption}\label{Assumption6}
There exists a $\nu\ge
0$ such that
%
%
\begin{equation} \label{201}
0<\inf_{|x|\ge1}\frac{|m(x)|}{|x|^\nu} \le
C \sup_{x\in R}\frac{|m(x)|}{1+|x|^\nu}<\infty,
\end{equation}
and there exist $0<\gamma'
\leq1$ and $\max\{0,3/4-2\nu\}<\gamma\le1$ such that
%
%
\begin{equation}\label{2011}\qquad
| m(x+y)-m(x)| \leq C|y|^{\gamma}\cases{
1+|x|^{\nu-1}+|y|^{\nu}, &\quad if $\nu>0$,\cr
1+|x|^{\gamma' -1}, &\quad if $\nu=0$,}
\end{equation}
for any $x,y\in R$ and for some constant $C>0$.
\end{Assumption}

Assumption \ref{Assumption6} is quite weak which is satisfied by a
large class of real functions
such as $m(x)=a_1+a_2x+\cdots+a_kx^{k-1}$,
$m(x) =a+b x^{c }$ and $m(x)=(a +b e^{x})/(1+e^{x})$. If $m(x)$ is
positive(or negative) on $R$,
condition (\ref{2011}) is not necessary.
%
%
\begin{theorem}
\label{th4} In addition to Assumptions \ref{Assumption1}--\ref
{Assumption6}, $\int(1+|x|^{2\nu
+2})K(x)\,dx<\infty$
and $E|\ep_0|^{4\nu+2}<\infty$. Then,
under $H_1$, we have
%
%
\begin{equation}\label{205}
\lim_{n\to\infty}
P\biggl(\frac{S_{n}}{\sqrt{2} V_{n}} \ge t_{\al}\biggr)=1
\end{equation}
for any $\rho_n$ satisfying
$\rho_n^2n^{1/2+\nu}h^{1/2}\to\infty$, and for any $0<\al<1$, where
$\Phi(t_\al)=1-\al$ and $\Phi$ is the standard normal c.d.f.
\end{theorem}

Theorem \ref{th4} shows that our test has nontrivial power against
the local alternative whenever $\rho_n\to0$ at a rate that is slower
than $n^{-1/8-\nu/2}$, as $nh^2\to\infty$.
This is different
from the stationary situation where in general a test has a
nontrivial power if only $\rho_n\to0$ at a rate that is slower than~$n^{-1/2}$. It is interesting to notice that the rate is related to
the magnitude of $m(x)$ and the bandwidth $h$.
The test has stronger discriminatory power the larger the value of
$v$. The reason is that the nonlinear shape characteristics in $m(x)$
are magnified over a wide domain and this property is exploited by the
test because the nonstationary regressor is recurrent.

Theorem \ref{th4} seems to be new to the literature. Under very strict
restrictions (namely that $%
x_{t}$ is a random walk and $x_{t}$ is independent of $u_{t}$), the
result in Theorem \ref{th3} has been considered in
\citet{Gaoetal09N1}. Not only the generalization of our result,
but the techniques used in this paper are quite different from Gao
et~al. (\citeyear{Gaoetal09N1}, \citeyear{Gaoetal09N2}). To outline the
essentials of the argument in the proof of Theorem~\ref{th3}, under the
null hypothesis, we split $S_{n}$ as
%
%
\begin{eqnarray} \label{s1} S_n
&=& 2 \sum_{t=2}^{n}u_{t+1} Y_{nt} +2 \mathop{\sum_{i, t=1}}_{
i\not=t}^nu_{i+1}
[f(x_t, \theta)-f(x_t, \hat\theta)]K[(x_{t}-x_i)/h]\nonumber\\
&&{} + \mathop{\sum_{i, t=1}}_{i\not=t}^n
[f(x_i, \theta)-f(x_i, \hat\theta)]
[f(x_t, \theta)-f(x_t, \hat\theta)]K[(x_{t}-x_i)/h] \\
&=& 2 S_{1n}+2S_{2n}+S_{3n} \qquad\mbox{say},\nonumber
\end{eqnarray}
where $Y_{nt}=\sum_{i=1}^{t-1}u_{i+1}K[(x_{t}-x_{i})/h]$.
It will be proved in Section \ref{sec61} that terms $S_{2n}$ and $S_{3n}$ are
negligible in comparison with $S_{1n}$.
Furthermore it will be proved that, under the null hypothesis,
%
%
\begin{eqnarray} \label{lp1}
V_{n}^{2} &=& \sigma^{4} \mathop{\sum_{t,s=1}}_{t\not
=s}^{n}K^{2}[%
(x_{t}-x_{s})/h]+o_{P}(n^{3/2} h) \nonumber\\[-8pt]\\[-8pt]
&=& 2\sigma^{2} \sum_{t=2}^{n}Y_{nt}^{2}+o_{P}(n^{3/2} h).\nonumber
\end{eqnarray}
By virtue of these facts, Theorem \ref{th3} follows from the following theorem,
giving a joint convergence result for $S_{1n}$ and its conditional
variance $\sum_{t=2}^{n}Y_{nt}^{2}$. This result, along with the
following Corollary \ref{cor1},
is of some independent interest.
%
%
\begin{theorem}
\label{th5} Under Assumptions \ref{Assumption1}--\ref{Assumption3},
$nh^{2}\rightarrow\infty$ and $%
nh^{4}\log^{2}n\rightarrow0$, we have
%
%
\begin{equation} \label{192}
\Biggl(\frac{1}{\sigma d_{n}}\sum_{t=2}^{n}u_{t+1} Y_{nt}, \frac
{1}{%
d_{n}^{2}}\sum_{t=2}^{n}Y_{nt}^{2}\Biggr)\rightarrow_{D}(\eta
N,\eta^{2}),
\end{equation}
where $d_{n}^{2}=(2 \phi)^{-1}\sigma^{2}n^{3/2}h \int_{-\infty
}^{\infty}K^{2}(x)\,dx$, $\eta
^{2}=L_{G}(1,0)$ is the self intersection local time generated by the
process $G=\int_{0}^{t}e^{\kappa(t-s)}\,dW(s)$, and
$N$ is a~standard normal variate which is independent of $\eta
^{2}$.
\end{theorem}
%
%
\begin{cor} \label{cor1} Under Assumptions \ref{Assumption1}--\ref
{Assumption5}, we have
\[
\frac{S_{n}}{\tau_{n}}\rightarrow_{D}\eta N,
\]
where $\tau_{n}^{2}=(8\phi)^{-1} \sigma
^{4} n^{3/2}h \int_{-\infty}^{\infty}K^{2}(x)\,dx$, $\eta
^{2}$ and $N$ are defined as in Theorem~\ref{th5}.
\end{cor}

Here and below, we define
%
%
\begin{eqnarray} \label{tt2}
L_{G}(t,u) &=&\lim_{\varepsilon\rightarrow0}\frac{1}{2\varepsilon}%
\int_{0}^{t}\int_{0}^{t}\mathbf{1}\bigl[\bigl\vert\bigl(
G(x)-G(y)\bigr)
-u\bigr\vert<\varepsilon\bigr]\,dx\,dy \nonumber\\[-8pt]\\[-8pt]
&=& \int_{0}^{t}\int_{0}^{t}\delta_{u}[G(x)-G(y)]\,dx\,dy,\nonumber
\end{eqnarray}
where $\delta_{u}$ is the dirac function. $L_{G}(t,u)$
characterizes the
amount of time over the interval $[ 0,t] $ that the
process $%
G( t) $ spends at a distance $u$ from itself, and is well
defined, as shown in Section \ref{sec5}. When $u=0$, $%
L_{G}(t,0)$ describes the self-intersection time of the process
$G( t)$. Using the definition of the dirac function, the
extended occupation times formula [e.g., \citet{RevYor99}, page 232],
and integration by parts with the local time measure, we may write
%
%
\begin{eqnarray}\label{tt3}
L_{G}(t,0) &=&2\int_{0}^{t}\int_{0}^{y}\delta_{0}[G(x)-G(y)]\,dx\,dy
\nonumber\\
&=&2\int_{0}^{t}\ell_{G}( s,G( s) ) \,ds
\nonumber\\[-8pt]\\[-8pt]
&=&2\int_{-\infty}^{\infty}\int_{0}^{t}\ell_{G}( s,a)\,
d\ell_{G}
( s,a) \,da \nonumber\\
&=&\int_{-\infty}^{\infty}\ell_{G}( t,a) ^{2}\,da,
\nonumber
\end{eqnarray}
where $\ell_{G}( t,a) $ is the local time spent by the
process $%
G $ at $a$ over the time interval $[ 0,t]$, namely,%
\[
\ell_{G}( t,a) =\int_{0}^{t}\delta
_{a}[G(s)]\,ds=\lim_{\varepsilon\rightarrow0}\frac{1}{2\varepsilon}%
\int_{0}^{t}\mathbf{1}[\vert G(s)-a\vert<\varepsilon
]\,ds.
\]
The process $\ell_{G}( s,G( s) ) $ is the
local time that the process $G$ has spent at its current position
$G( s) $ over the time interval $[ 0,s]$. It
appears in the limit theory for nonparametric nonstationary spurious
regression [\citet{Phi09}]. \citet{Ald86} gave (\ref{tt3}) for the
case of Brownian motion.

It is interesting to note that $S_{1n}$ is a martingale sequence with
conditional variance $\sum_{t=2}^{n}Y_{nt}^{2}$, suggesting that some
version of the martingale central limit theorem [e.g.,
\citet{HalHey80}, Chapter 3] may be applicable. However, the
problem is complicated by the U-statistic structure and the weak
convergence of the conditional variance, and use of existing limit
theory seems difficult. To investigate the asymtotics of $S_{1n}$, we
therefore develop our own approach. As part of this development, in
Section \ref{sec5}, we provide a general weak convergence theory to
intersection local time, which is of independent interest and useful in
other applications. The conditions required for this development are
weaker than those in establishing Theorem \ref{th5} and that section
may be read separately.

We finally remark that the restrictive condition on the bandwidth
$h$ in Theorems \ref{th3}--\ref{th5} (i.e., $nh^4\log^2 n\to0$) is mainly
used to offset the impact of the error terms in (3.5) and (3.6). It
seems difficult to relax this condition under the prevailing
Assumption \ref{Assumption2}, which allows for endogenity in the
regressor~$x_t$.
See, for instance, the proof of Proposition \ref{prop6} given
in the supplementary material [\citet{WanPhi}]. The restriction $nh^4\log^2 n\to0$ on~$h$
in Theorems \ref{th3}--\ref{th5}, however, can be reduced to the minimal
requirement $h\to0$, if Assumption \ref{Assumption2} is replaced by the
following
Assumption \ref{Assumption2zv}.

{\renewcommand{\theAssumption}{$2^*$}
%
%
\begin{Assumption}\label{Assumption2zv}
For\vspace*{1pt} each $n\ge1$,
$\{u_{t},\mathcal{F}_{t, n}\}_{1\le t\le n}$ forms a
martingale difference satisfying ${\lim_{t\to\infty}\sup_{n\ge
t}}|E(u_{t+1}^{2}\mid%
\mathcal{F}_{t, n})-\sigma^{2}| =0$, a.s. and
\[
\sup
_{n\ge t\geq1}E(|u_{t+1}|^{4}\mid\mathcal{F}%
_{t, n})<\infty,
\]
where
\[
\mathcal{F}_{t, n} =\si(u_1,\ldots,u_t; x_1,\ldots,x_n),\qquad t=1,2,\ldots,n;
n\ge1.
\]
\end{Assumption}}

Note that Assumption \ref{Assumption2zv} holds true if $x_t$ is
independent of $u_t$,
and $\{u_{t},\allowbreak \mathcal{F}%
_{t}\}_{t\geq1}$ forms a martingale difference satisfying
$E(u_{t+1}^{2}\mid%
\mathcal{F}_{t})\rightarrow_{\mathrm{a.s.}}\sigma^{2}>0$ as $t\rightarrow
\infty$ and $\sup_{t\geq1}E(|u_{t+1}|^{4}\mid\mathcal{F}
_{t})<\infty$, where $\mathcal{F}_{t}$ is a sequence of increasing
$\sigma$-fields. The independence assumption was used in
Gao et~al. (\citeyear{Gaoetal09N1}) to establish a similar version of Theorem \ref{th3}.

\section{Simulations}\label{sec4}

Simulations were conducted to evaluate the finite sample performance of the
statistic $Z_{n}=S_{n}/\sqrt{2}V_{n}$ under the null and some local
alternatives under various assumptions about the generating mechanism.
The results are summarized here, and more detailed
findings are reported in the supplementary material [\citet{WanPhi}]. The model
followed~(2.1) with $%
y_{t+1}=f( x_{t}) +u_{t+1}$, $x_{t}=x_{t-1}+\eta_{t}$, $x_{0}=0$,
and~$\eta_{t}$ generated by an AR(1) process $\eta_{t}=\lambda\eta
_{t-1}+\varepsilon_{t}$ or an MA(1) process $\eta_{t}=\varepsilon
_{t}+\lambda\varepsilon_{t-1}$ with $(
u_{t},\varepsilon_{t}) \sim\mbox{i.i.d. }N\bigl( 0,\bigl(
{1 \atop r} \enskip
{r \atop 1}
\bigr) \bigr)$. A linear null hypothesis $H_{0}\dvtx f(
x)
=\theta_{0}+\theta_{1}x$ was used together with polynomial local
alternatives $H_{1}\dvtx f( x) =\theta_{0}+\theta
_{1}x+\rho
_{n}\vert x\vert^{\nu}$, with $\rho_{n}=1/(n^{
{1}/{4}+%
{\nu}/{3}}h^{{1}/{4}})$. The parameter settings were $\theta
_{0}=0,\theta_{1}=1$, $\nu\in\{ 0.5,1.5,2,3\}$ and
$r\in\{
0,\pm0.5,\pm0.75\}$. Results are reported for sample sizes
$n\in
\{ 100,200,500\} $ and bandwidth settings $h=n^{-p}$ for
$p\in
\{ \frac{1}{4},\frac{1}{3},\frac{1}{2.5}\} $. Note that
$%
h=n^{-1/4}$ satisfies Assumption \ref{Assumption2zv} but not
Assumption~\ref{Assumption2}. The number of replications was $5000$.

%
\begin{table}[t]
\tablewidth=300pt
\caption{Size: $\eta_{t}=\varepsilon_{t}$}\label{table1}
\begin{tabular*}{\tablewidth}{@{\extracolsep{\fill}}
llccccccc@{}}
\hline
&& \multicolumn{3}{c}{\textbf{Nominal size} $\bolds{5\%}$} & &
\multicolumn{3}{c@{}}{\textbf{Nominal size} $\bolds{1\%}$} \\[-4pt]
&& \multicolumn{3}{c}{\hrulefill} & &
\multicolumn{3}{c@{}}{\hrulefill} \\
$\bolds{n}$ & \multicolumn{1}{l}{$\bolds{h=}${\hspace*{-9pt}}} &
\multicolumn{1}{c}{$\bolds{n^{-{1}/{4}}}$} &
\multicolumn{1}{c}{$\bolds{n^{-{1}/{3}}}$} &
\multicolumn{1}{c}{$\bolds{n^{-{1}/{2.5}}}$} & &
\multicolumn{1}{c}{$\bolds{n^{-{1}/{4}}}$} &
\multicolumn{1}{c}{$\bolds{n^{-{1}/{3}}}$} &
\multicolumn{1}{c@{}}{$\bolds{n^{-{1}/{2.5}}}$} \\
\hline
\multicolumn{9}{@{}c@{}}{$r=0$}\\[3pt]
100 & & 0.028 & 0.035 & 0.033 & & 0.006 & 0.006 & 0.007 \\
200 & & 0.034 & 0.042 & 0.041 & & 0.007 & 0.007 & 0.008 \\
500 & & 0.044 & 0.045 & 0.050 & & 0.009 & 0.010 & 0.010 \\ 
[3pt]
\multicolumn{9}{@{}c@{}}{$r=0.5$}\\ 
[3pt]
100 & & 0.030 & 0.035 & 0.040 & & 0.006 & 0.007 & 0.007 \\
200 & & 0.038 & 0.044 & 0.045 & & 0.009 & 0.008 & 0.008 \\
500 & & 0.041 & 0.045 & 0.048 & & 0.008 & 0.009 & 0.009 \\ 
[3pt]
\multicolumn{9}{@{}c@{}}{$r=-0.5$}\\ 
[3pt]
100 & & 0.031 & 0.035 & 0.037 & & 0.007 & 0.008 & 0.008 \\
200 & & 0.036 & 0.045 & 0.046 & & 0.007 & 0.008 & 0.009 \\
500 & & 0.041 & 0.047 & 0.051 & & 0.009 & 0.010 & 0.011 \\
\hline
\end{tabular*}
\end{table}

%
\begin{table}
\tablewidth=300pt
\caption{Size: $\eta_{t}=\lambda\eta_{t-1}+\varepsilon_{t}$, $r=\pm0.5$}\label{table2}
\begin{tabular*}{\tablewidth}{@{\extracolsep{\fill}}
llccccccc@{}}
\hline
&& \multicolumn{3}{c}{\textbf{Nominal size} $\bolds{5\%}$} & &
\multicolumn{3}{c@{}}{\textbf{Nominal size} $\bolds{1\%}$} \\[-4pt]
&& \multicolumn{3}{c}{\hrulefill} & &
\multicolumn{3}{c@{}}{\hrulefill} \\
$\bolds{n}$ & \multicolumn{1}{l}{$\bolds{h=}${\hspace*{-9pt}}} &
\multicolumn{1}{c}{$\bolds{n^{-{1}/{4}}}$} &
\multicolumn{1}{c}{$\bolds{n^{-{1}/{3}}}$} &
\multicolumn{1}{c}{$\bolds{n^{-{1}/{2.5}}}$} & &
\multicolumn{1}{c}{$\bolds{n^{-{1}/{4}}}$} &
\multicolumn{1}{c}{$\bolds{n^{-{1}/{3}}}$} &
\multicolumn{1}{c@{}}{$\bolds{n^{-{1}/{2.5}}}$} \\
\hline
\multicolumn{9}{@{}c@{}}{$r=0.5$, $\lambda=0.4$}\\
[3pt]
100 & & 0.034 & 0.038 & 0.041 & & 0.002 & 0.004 & 0.005 \\
200 & & 0.044 & 0.044 & 0.047 & & 0.004 & 0.006 & 0.007 \\
500 & & 0.058 & 0.058 & 0.057 & & 0.007 & 0.010 & 0.011 \\ 
[3pt]
\multicolumn{9}{@{}c@{}}{$r=0.5$, $\lambda=-0.4$}\\
[3pt]
100 & & 0.038 & 0.042 & 0.046 & & 0.013 & 0.013 & 0.011 \\
200 & & 0.051 & 0.051 & 0.051 & & 0.018 & 0.015 & 0.014 \\
500 & & 0.070 & 0.061 & 0.057 & & 0.026 & 0.022 & 0.016 \\ 
[3pt]
\multicolumn{9}{@{}c@{}}{$r=-0.5$, $\lambda=0.4$}\\ 
[3pt]
100 & & 0.034 & 0.038 & 0.040 & & 0.002 & 0.004 & 0.005 \\
200 & & 0.044 & 0.044 & 0.048 & & 0.004 & 0.006 & 0.007 \\
500 & & 0.058 & 0.058 & 0.057 & & 0.007 & 0.009 & 0.011 \\ 
[3pt]
\multicolumn{9}{@{}c@{}}{$r=-0.5$, $\lambda=-0.4$}\\ 
[3pt]
100 & & 0.035 & 0.040 & 0.043 & & 0.012 & 0.012 & 0.012 \\
200 & & 0.050 & 0.049 & 0.050 & & 0.018 & 0.015 & 0.013 \\
500 & & 0.073 & 0.064 & 0.056 & & 0.026 & 0.018 & 0.016 \\
\hline
\end{tabular*}
\end{table}

Table \ref{table1} shows the actual size of the test for various $n$
and bandwidth
choices $h$ and for both exogenous $( r=0) $ and
endogenous $%
( r=\pm0.5) $ regressor cases with serially uncorrelated
errors ($\lambda=0$). Table \ref{table2} shows the corresponding
results for AR errors
with $\lambda=\pm0.4$. Size results for MA errors are similar and are
given in the supplementary material [\citet{WanPhi}]. Under i.i.d. errors the test is somewhat undersized
for $n=100,200$ but is close to the nominal for $n=500$ and for all
bandwidth choices. There is some mild oversizing under serially
dependent $%
\eta_{t}$ when $\lambda=-0.4$ for bandwidth $h=n^{-1/4}$, but size seems
satisfactory for $\lambda=0.4$ and for the smaller bandwidths $%
h=n^{-1/3},n^{-1/2.5}$. Since negative $\lambda$ reduces the long run
moving average coefficient $\phi$ [$\phi=1/( 1-\lambda)
$ for
AR $\eta_{t}$] these results suggest that the strength of the long run
signal in $x_{t}$ (measured by the long-run variance of~$\eta_{t}$) affects
the performance of the test. On the other hand, endogeneity at the
correlation level $r=\pm0.5$ appears to have little effect on performance,
which mirrors results for estimation in the nonlinear nonstationary case
[\citet{WanPhi09N2}]. Higher levels of correlation ($r=\pm0.75$)
produce some size distortion when there is serial dependence, but not
when the errors are independent; see Table \ref{table3}.

Table \ref{table4}--\ref{table6} show test power against the local
alternative $H_{1}$ for polynomial alternatives (cubic $\nu=3$,
quadratic $\nu=2$ and three halves $\nu=1.5$). Results for the case
$\nu=0.5$ are given in the supplementary material
[\citet{WanPhi}]. Again, there is little difference between the
exogenous and endogenous cases, so only the endogenous case is reported
here. As may be expected, there is greater local discriminatory power
for cubic
($\nu=3$%
) than quadratic ($\nu=2$) or three halves ($\nu=1.5$) alternatives.
For $%
n=100$ (500) power is greater than $69\%$ ($90\%$) for a nominal
$1\%$
test and greater than $74\%$ (92\%) for a nominal $5\%$ test when $\nu=3$
under AR errors with $\lambda=0.4$ (Table \ref{table4}). The
corresponding results
when $\nu=2$ and $n=100$ (500) are 15\% (38\%) for a nominal 1\% test and
23\% (46\%) for a nominal 5\% test (Table \ref{table5}). Serial
dependence affects
power, which is higher for $\lambda=0.4$ than for $\lambda=-0.4$ in all
cases. So lower long-run signal strength in the regressor tends to reduce
discriminatory power. For $\nu=1.5$ and $\lambda=-0.4$, power is low even
for $n=500$ (2$+$\% for a 1\% test and 7$+$\% for a 5\% test, Table \ref
{table6}). Low
power also occurs against the local alternative with $\nu=0.5$ [see
\citet{WanPhi}], which also reduces signal strength in the regressor
function. Thus,
discriminatory power is dependent on the specific alternative and, as
asymptotic theory suggests, is sensitive to the magnitude rate ($\nu$)
of $%
m( x) $ as $\vert x\vert\rightarrow\infty$.

%
\begin{table}
\tablewidth=300pt
\caption{Size: $\eta_{t}=\lambda\eta_{t-1}+\varepsilon_{t}$, $r=\pm0.75$}\label{table3}
\begin{tabular*}{\tablewidth}{@{\extracolsep{\fill}}
llccccccc@{}}
\hline
&& \multicolumn{3}{c}{\textbf{Nominal size} $\bolds{5\%}$} & &
\multicolumn{3}{c@{}}{\textbf{Nominal size} $\bolds{1\%}$} \\[-4pt]
&& \multicolumn{3}{c}{\hrulefill} & &
\multicolumn{3}{c@{}}{\hrulefill} \\
$\bolds{n}$ & \multicolumn{1}{l}{$\bolds{h=}${\hspace*{-9pt}}} &
\multicolumn{1}{c}{$\bolds{n^{-{1}/{4}}}$} &
\multicolumn{1}{c}{$\bolds{n^{-{1}/{3}}}$} &
\multicolumn{1}{c}{$\bolds{n^{-{1}/{2.5}}}$} & &
\multicolumn{1}{c}{$\bolds{n^{-{1}/{4}}}$} &
\multicolumn{1}{c}{$\bolds{n^{-{1}/{3}}}$} &
\multicolumn{1}{c@{}}{$\bolds{n^{-{1}/{2.5}}}$} \\
\hline
\multicolumn{9}{@{}c@{}}{$r=0.75$, $\lambda=0.4$}\\
[3pt]
100 & & 0.036 & 0.038 & 0.039 & & 0.003 & 0.003 & 0.004 \\
200 & & 0.043 & 0.049 & 0.050 & & 0.005 & 0.006 & 0.007 \\
500 & & 0.057 & 0.055 & 0.053 & & 0.007 & 0.009 & 0.008 \\ 
[3pt]
\multicolumn{9}{@{}c@{}}{$r=0.75$, $\lambda=-0.4$}\\
[3pt]
100 & & 0.074 & 0.068 & 0.027 & & 0.036 & 0.033 & 0.027 \\
200 & & 0.108 & 0.096 & 0.087 & & 0.050 & 0.043 & 0.034 \\
500 & & 0.177 & 0.140 & 0.115 & & 0.094 & 0.062 & 0.048 \\ 
[3pt]
\multicolumn{9}{@{}c@{}}{$r=0.75$, $\lambda=0$} \\
[3pt]
100 & & 0.026 & 0.029 & 0.032 & & 0.005 & 0.006 & 0.006 \\
200 & & 0.037 & 0.044 & 0.046 & & 0.007 & 0.008 & 0.010 \\
500 & & 0.040 & 0.042 & 0.047 & & 0.008 & 0.009 & 0.009 \\ 
[3pt]
\multicolumn{9}{@{}c@{}}{$r=-0.75$, $\lambda=0$} \\
[3pt]
100 & & 0.027 & 0.035 & 0.036 & & 0.005 & 0.008 & 0.007 \\
200 & & 0.036 & 0.040 & 0.043 & & 0.008 & 0.010 & 0.010 \\
500 & & 0.041 & 0.045 & 0.044 & & 0.008 & 0.008 & 0.009 \\ 
[3pt]
\multicolumn{9}{@{}c@{}}{$r=-0.75$, $\lambda=0.4$} \\
[3pt]
100 & & 0.074 & 0.071 & 0.063 & & 0.003 & 0.004 & 0.004 \\
200 & & 0.103 & 0.085 & 0.074 & & 0.011 & 0.012 & 0.011 \\
500 & & 0.135 & 0.105 & 0.088 & & 0.027 & 0.020 & 0.015 \\ 
[3pt]
\multicolumn{9}{@{}c@{}}{$r=-0.75$, $\lambda=-0.4$} \\
[3pt]
100 & & 0.070 & 0.066 & 0.065 & & 0.033 & 0.026 & 0.023 \\
200 & & 0.109 & 0.094 & 0.087 & & 0.055 & 0.042 & 0.033 \\
500 & & 0.175 & 0.136 & 0.109 & & 0.093 & 0.065 & 0.048 \\
\hline
\end{tabular*}
\end{table}

Overall, the finite sample results reflect the asymptotic theory and seem
reasonable for practical use in testing when there is some endogeneity in
nonparametric nonstationary regression, especially if smaller bandwidth
choices than usual are employed. In cases of serial dependence when the
long-run signal strength in the regressor $x_{t}$ is reduced, finite sample
adjustments for the test critical values may be useful in correcting size,
as has been found for i.i.d. and stationary regressors [\citet{LiWan98}].

%
\begin{table}
\tablewidth=300pt
\caption{Local power: $\nu=3$, $\eta _{t}=\lambda \eta_{t-1}+\varepsilon_{t}$,
$r=\pm0.5$}\label{table4}
\begin{tabular*}{\tablewidth}{@{\extracolsep{\fill}}
llccccccc@{}}
\hline
&& \multicolumn{3}{c}{\textbf{Nominal size} $\bolds{5\%}$} & &
\multicolumn{3}{c@{}}{\textbf{Nominal size} $\bolds{1\%}$} \\[-4pt]
&& \multicolumn{3}{c}{\hrulefill} & &
\multicolumn{3}{c@{}}{\hrulefill} \\
$\bolds{n}$ & \multicolumn{1}{l}{$\bolds{h=}${\hspace*{-9pt}}} &
\multicolumn{1}{c}{$\bolds{n^{-{1}/{4}}}$} &
\multicolumn{1}{c}{$\bolds{n^{-{1}/{3}}}$} &
\multicolumn{1}{c}{$\bolds{n^{-{1}/{2.5}}}$} & &
\multicolumn{1}{c}{$\bolds{n^{-{1}/{4}}}$} &
\multicolumn{1}{c}{$\bolds{n^{-{1}/{3}}}$} &
\multicolumn{1}{c@{}}{$\bolds{n^{-{1}/{2.5}}}$} \\
\hline
\multicolumn{9}{@{}c@{}}{$r=0.5$, $\lambda=0.4$}\\
[3pt]
100 & & 0.819 & 0.779 & 0.743 & & 0.787 & 0.739 & 0.693 \\
200 & & 0.906 & 0.878 & 0.845 & & 0.892 & 0.849 & 0.811 \\
500 & & 0.971 & 0.950 & 0.923 & & 0.963 & 0.935 & 0.901 \\ 
[3pt]
\multicolumn{9}{@{}c@{}}{$r=0.5$, $\lambda=-0.4$} \\
[3pt]
100 & & 0.247 & 0.211 & 0.179 & & 0.197 & 0.154 & 0.126 \\
200 & & 0.358 & 0.306 & 0.265 & & 0.302 & 0.247 & 0.199 \\
500 & & 0.522 & 0.448 & 0.389 & & 0.458 & 0.376 & 0.310 \\ 
[3pt]
\multicolumn{9}{@{}c@{}}{$r=-0.5$, $\lambda=0.4$} \\
[3pt]
100 & & 0.829 & 0.780 & 0.743 & & 0.792 & 0.742 & 0.696 \\
200 & & 0.910 & 0.879 & 0.845 & & 0.891 & 0.851 & 0.813 \\
500 & & 0.965 & 0.947 & 0.921 & & 0.957 & 0.931 & 0.903 \\ 
[3pt]
\multicolumn{9}{@{}c@{}}{$r=-0.5$, $\lambda=-0.4$} \\
[3pt]
100 & & 0.238 & 0.204 & 0.176 & & 0.189 & 0.151 & 0.127 \\
200 & & 0.352 & 0.297 & 0.253 & & 0.295 & 0.239 & 0.193 \\
500 & & 0.513 & 0.431 & 0.367 & & 0.449 & 0.367 & 0.301 \\
\hline
\end{tabular*}
\end{table}

%
\begin{table}
\tablewidth=300pt
\caption{Local power: $\nu=2$, $\eta _{t}=\lambda \eta_{t-1}+\varepsilon_{t}$,
$r=\pm0.5$}\label{table5}
\begin{tabular*}{\tablewidth}{@{\extracolsep{\fill}}
llccccccc@{}}
\hline
&& \multicolumn{3}{c}{\textbf{Nominal size} $\bolds{5\%}$} & &
\multicolumn{3}{c@{}}{\textbf{Nominal size} $\bolds{1\%}$} \\[-4pt]
&& \multicolumn{3}{c}{\hrulefill} & &
\multicolumn{3}{c@{}}{\hrulefill} \\
$\bolds{n}$ & \multicolumn{1}{l}{$\bolds{h=}${\hspace*{-9pt}}} &
\multicolumn{1}{c}{$\bolds{n^{-{1}/{4}}}$} &
\multicolumn{1}{c}{$\bolds{n^{-{1}/{3}}}$} &
\multicolumn{1}{c}{$\bolds{n^{-{1}/{2.5}}}$} & &
\multicolumn{1}{c}{$\bolds{n^{-{1}/{4}}}$} &
\multicolumn{1}{c}{$\bolds{n^{-{1}/{3}}}$} &
\multicolumn{1}{c@{}}{$\bolds{n^{-{1}/{2.5}}}$} \\
\hline
\multicolumn{9}{@{}c@{}}{$r=0.5$, $\lambda=0.4$}\\
[3pt]
100 & & 0.357 & 0.282 & 0.228 & & 0.282 & 0.205 & 0.147 \\
200 & & 0.484 & 0.389 & 0.315 & & 0.418 & 0.310 & 0.228 \\
500 & & 0.682 & 0.557 & 0.458 & & 0.616 & 0.482 & 0.376 \\ 
[3pt]
\multicolumn{9}{@{}c@{}}{$r=0.5$, $\lambda=-0.4$}\\
[3pt]
100 & & 0.058 & 0.054 & 0.053 & & 0.027 & 0.020 & 0.016 \\
200 & & 0.103 & 0.083 & 0.068 & & 0.048 & 0.034 & 0.024 \\
500 & & 0.169 & 0.118 & 0.094 & & 0.098 & 0.057 & 0.036 \\ 
[3pt]
\multicolumn{9}{@{}c@{}}{$r=-0.5$, $\lambda=0.4$} \\
[3pt]
100 & & 0.114 & 0.123 & 0.128 & & 0.065 & 0.066 & 0.067 \\
200 & & 0.226 & 0.235 & 0.244 & & 0.157 & 0.159 & 0.160 \\
500 & & 0.437 & 0.457 & 0.462 & & 0.350 & 0.359 & 0.367 \\ 
[3pt]
\multicolumn{9}{@{}c@{}}{$r=-0.5$, $\lambda=-0.4$} \\
[3pt]
100 & & 0.056 & 0.050 & 0.046 & & 0.022 & 0.016 & 0.014 \\
200 & & 0.102 & 0.082 & 0.066 & & 0.053 & 0.031 & 0.022 \\
500 & & 0.173 & 0.123 & 0.096 & & 0.103 & 0.061 & 0.037 \\
\hline
\end{tabular*}
\end{table}

In practice, the exact $\alpha$-level critical value $\ell_{\alpha
}(
h) $ ($0<\alpha<1$) of the finite sample distribution of
$S_{n}/\sqrt{%
2}V_{n}$ depends on all the unknown parameters and functions in the
model. The development of a rigorous theory of approximation for $\ell
_{\alpha }( h) $ and the choice of an optimal bandwidth for use in
testing are challenging problems in the nonstationary setting.
\citet{Gaoetal09N1} provided an approximate value of
$\ell_{\alpha}( h) $ by using the bootstrap and considered numerical
solutions for a bandwidth $h$ that optimizes the power function, both
under the assumption that $x_{t}$
and $%
u_{t}$ are independent. It is not clear at the moment whether similar
techniques can be rigorously justified in the current general model and
there is presently no optimal approach to bandwidth selection. The
%
%
\begin{table}
\tablewidth=300pt
\caption{Local power: $\nu=1.5$, $\eta _{t}=\lambda \eta_{t-1}+\varepsilon_{t}$,
$r=\pm0.5$}\label{table6}
\begin{tabular*}{\tablewidth}{@{\extracolsep{\fill}}
llccccccc@{}}
\hline
&& \multicolumn{3}{c}{\textbf{Nominal size} $\bolds{5\%}$} & &
\multicolumn{3}{c@{}}{\textbf{Nominal size} $\bolds{1\%}$} \\[-4pt]
&& \multicolumn{3}{c}{\hrulefill} & &
\multicolumn{3}{c@{}}{\hrulefill} \\
$\bolds{n}$ & \multicolumn{1}{l}{$\bolds{h=}${\hspace*{-9pt}}} &
\multicolumn{1}{c}{$\bolds{n^{-{1}/{4}}}$} &
\multicolumn{1}{c}{$\bolds{n^{-{1}/{3}}}$} &
\multicolumn{1}{c}{$\bolds{n^{-{1}/{2.5}}}$} & &
\multicolumn{1}{c}{$\bolds{n^{-{1}/{4}}}$} &
\multicolumn{1}{c}{$\bolds{n^{-{1}/{3}}}$} &
\multicolumn{1}{c@{}}{$\bolds{n^{-{1}/{2.5}}}$} \\
\hline
\multicolumn{9}{@{}c@{}}{$r=0.5$, $\lambda=0.4$}\\
[3pt]
100 & & 0.058 & 0.051 & 0.045 & & 0.021 & 0.012 & 0.010 \\
200 & & 0.087 & 0.065 & 0.057 & & 0.040 & 0.022 & 0.015 \\
500 & & 0.158 & 0.103 & 0.077 & & 0.096 & 0.046 & 0.024 \\ 
[3pt]
\multicolumn{9}{@{}c@{}}{$r=0.5$, $\lambda=-0.4$}\\
[3pt]
100 & & 0.043 & 0.040 & 0.041 & & 0.016 & 0.014 & 0.012 \\
200 & & 0.061 & 0.058 & 0.055 & & 0.024 & 0.019 & 0.015 \\
500 & & 0.096 & 0.074 & 0.070 & & 0.038 & 0.031 & 0.023 \\ 
[3pt]
\multicolumn{9}{@{}c@{}}{$r=-0.5$, $\lambda=0.4$} \\
[3pt]
100 & & 0.066 & 0.053 & 0.050 & & 0.025 & 0.015 & 0.011 \\
200 & & 0.093 & 0.065 & 0.052 & & 0.046 & 0.023 & 0.015 \\
500 & & 0.152 & 0.094 & 0.090 & & 0.088 & 0.042 & 0.023 \\ 
[3pt]
\multicolumn{9}{@{}c@{}}{$r=-0.5$, $\lambda=-0.4$} \\
[3pt]
100 & & 0.049 & 0.049 & 0.049 & & 0.018 & 0.017 & 0.013 \\
200 & & 0.063 & 0.058 & 0.059 & & 0.024 & 0.021 & 0.017 \\
500 & & 0.092 & 0.074 & 0.064 & & 0.037 & 0.029 & 0.021 \\
\hline
\end{tabular*}
\vspace*{3pt}
\end{table}
investigation of such finite sample adjustments and selection criteria
is therefore left for later research. Earlier analysis of the
restrictions on the bandwidth in Theorems \ref{th3}--\ref{th5}, in conjunction with
the simulation evidence, indicates that smaller bandwidths than usual
for stationary regression are likely to be more reliable in practical
work for specification testing of nonlinear nonstationary regression.%

\section{Convergence to intersection local time}\label{sec5}

Consider a linear process $\{\eta_{j}$, $j\geq1\}$ defined by $\eta
_{j}=\sum_{k=0}^{\infty} \phi_{k} \epsilon_{j-k}$, where $\{
\epsilon
_{j},j\in Z\}$ is a sequence of i.i.d. random variables with $E\epsilon_{0}=0$
and $E\epsilon_{0}^{2}=1$, and the coefficients $\phi_{k},k\geq0$ are
assumed to satisfy ${\sum_{k=0}^{\infty}}|\phi_{k}|<\infty$ and $\phi
\equiv\sum_{k=0}^{\infty}\phi_{k}\not=0$. Let
%
%
\begin{equation}\label{mod1}
y_{k,n}=\rho y_{k-1,n}+\eta_{k},\qquad y_{0,n}=0,\qquad \rho=1+\kappa/n,
\end{equation}
where $\kappa$ is a constant. The array $y_{k,n}$, $k\geq0$ is known
as a
nearly unstable process or, in the econometric literature, as a
near-integrated time series. Write $x_{k,n}=y_{k,n}/\sqrt{n}\phi$. The
classical invariance principle gives
%
%
\begin{equation}\label{no1}\qquad
x_{[nt],n}\Rightarrow G(t):=\int_{0}^{t}e^{\kappa
(t-s)}\,dW(s)=W(t)+\kappa
\int_{0}^{t}e^{\kappa(t-s)}W(s)\,ds
\end{equation}
on $D[0,1]$, where $W(t)$ is a standard Brownian motion [e.g.,
\citet{Phi87}, \citet{BucCha07}, \citet{WanPhi09N2}].
Furthermore, $%
\{\epsilon_{j},j\in Z\}$ can be redefined on a richer probability space
which also contains a standard Brownian motion $W_{1}(t)$ such that
%
%
\begin{equation} \label{no11}
\sup_{0\leq t\leq1}\bigl|x_{[nt],n}-G_{1}(t)\bigr|=o_{P}(1),
\end{equation}
where $G_{1}(t)=W_{1}(t)+\kappa\int_{0}^{t}e^{\kappa(t-s)}W_{1}(s)\,ds$.
Indeed, by noting on the richer space that
%
%
\begin{equation}
\sup_{0\leq t\leq1}\Biggl|\frac{1}{\sqrt{n}}\sum
_{j=1}^{[nt]}\epsilon
_{j}-W_{1}(t)\Biggr|=o_{P}(1)
\end{equation}
[see, e.g., \citet{CsoRev81}], and using this result
in place of the fact that $\frac{1}{\sqrt{n}}\sum
_{j=1}^{[nt]}\epsilon
_{j}\Rightarrow W(t)$ on $D[0,1]$, the same technique as in the proof of
\citet{Phi87} [see also \citet{ChaWei87}] yields
\[
\sup_{0\leq t\leq1}\Biggl|\frac{1}{\sqrt{n}}\sum_{j=1}^{[nt]}\rho
^{[
nt]-j}\epsilon_{j}-G_{1}(t)\Biggr|=o_{P}(1).
\]
The result (\ref{no11}) can now be obtained by the same argument, with
minor modifications, as in the proof of Proposition 7.1 in
\citet{WanPhi09N2}.

The aim of this section is to investigate the asymptotic behavior of a~functional $S_{[nr]}$ of the $x_{k,n}$, defined by
%
%
\begin{equation} \label{tt5}
S_{[nr]}=\frac{c_{n}}{n^{2}} \sum_{k,j=1}^{[nr]}g[c_{n}%
(x_{k,n}-x_{j,n})],
\end{equation}
where $g$ is a real function on $R$, and $c_{n}$ is a certain sequence of
positive constants. Under certain conditions on $g(x)$, $\epsilon_{0}$
and $%
c_{n}$, it is established that, for each fixed $0<r\le1$, $S_{[nr]}$ converges
to an intersection local time process of $G(t)$. Explicitly, we have the
following main result.
%
%
\begin{theorem}
\label{th1a} Suppose that $\int_{-\infty}^{\infty}|g(x)|\,dx<\infty
$, $%
\omega\equiv\int_{-\infty}^{\infty}g(x)\,dx\not=0$ and $\int
_{-\infty
}^{\infty}|Ee^{it\epsilon_{0}}|\,dt<\infty$. Then, for any $%
c_{n}\rightarrow\infty$, $n/c_{n}\rightarrow\infty$ and fixed $r\in
(0,1]$,
%
%
\begin{equation} \label{r1a7}
S_{[nr]}\rightarrow_{D}\omega L_{G}(r,0),
\end{equation}
where $L_{G}(t,u)$ is the intersection local time of $G(t)$ defined in
(\ref%
{tt2}). Furthermore, under the same probability space for which (%
\ref{no11}) holds, we have that, for any $c_{n}\rightarrow\infty$
and $%
n/c_{n}\rightarrow\infty$,
%
%
\begin{equation} \label{r28}
\sup_{0\leq r\leq1}\bigl|S_{[nr]}-\omega L_{G_{1}}(r,0)
\bigr|\rightarrow_{P}0.
\end{equation}
\end{theorem}

The integrability condition on the characteristic function of $\epsilon_{0}$
can be weakened if we place further restrictions on $g(x)$. Indeed, we have
the following theorem.
%
%
\begin{theorem}
\label{th2a} Theorem \ref{th1a} still holds if $\int_{-\infty
}^{\infty}|Ee^{it\epsilon_{0}}|\,dt<\infty$ is replaced by the
Cram\'er condition, that is, ${\limsup_{|t|\rightarrow\infty
}}|Ee^{it\epsilon_{0}}|<1$, and, in addition to the stated conditions
already on $g(x)$, we have $|g(x)|\leq M/(1+|x|^{1+b})$ for some
$b>0$, where $M$ is a constant.
\end{theorem}

It is interesting to notice that the additional condition on $g(x)$ in
Theorem \ref{th2a} cannot be reduced without further restriction on $%
\epsilon_{0}$ like that in Theorem \ref{th1a}. This claim can
be
explained as in Example 4.2.2 of \citet{BorIbr94} with
some minor modifications. On the other hand, the asymptotic behavior
of $S_{[nr]}$ when $c_{n}=1$ is quite different, as seen in the
following theorem.
%
%
\begin{theorem}
\label{th3a} Suppose that $g(x)$ is Borel measurable function satisfying
%
%
\begin{equation} \label{who19}
\lim_{h\rightarrow0}\int_{-K}^{K}{|x|^{\alpha-1}\sup_{|u|\leq
h}}|g(x+u)-g(x)|\,dx=0
\end{equation}
for all $K>0$ and some $0<\alpha\leq1$. Then, under the same probability
space for which (\ref{no11}) holds, we have
%
%
\begin{equation} \label{r1a7a}
\sup_{0\leq r\leq1}\Biggl|\frac{1}{n^{2}} \sum_{k,j=1}^{[nr]}g(
x_{k,n}-x_{j,n})-\int_{0}^{r}\int
_{0}^{r}g[G_{1}(u)-G_{1}(v)]\,du\,dv\Biggr|%
=o_{P}(1).\hspace*{-35pt}
\end{equation}
\end{theorem}

We mention that condition (\ref{who19}) is quite weak. Indeed, example
2.8 and the discussion following Theorem 2.3 in \citet{BerHor06}
shows that (\ref{who19}) cannot
be replaced by
\[
\lim_{h\rightarrow0}\int_{-K}^{K}|x|^{\alpha-1}|g(x+u)-g(x)|\,dx=0
\]
for all $K>0$ and some $0<\alpha\leq1$.

Local time has figured in much recent work on parametric and
nonparametric estimation with nonstationary data. Motivated by
nonlinear regression with integrated time series [Park and Phillips
(\citeyear{ParPhi99}, \citeyear{ParPhi01})] and nonparametric
estimation of nonlinear cointegration models, many authors
[\citet{PhiPar98}, \citet{KarTjs01},
\citet{KarMykTjs07},
\citet{WanPhi09N1}] have used or proved weak convergence to the
local time of a stochastic process, including results of the following
type: under certain conditions on the function $g$, the limiting
stochastic process $G(t)$, a sequence $c_{n}\rightarrow\infty$, and
normalized data $x_{k,n}$
%
%
\begin{equation} \label{aa0}
\frac{c_{n}}{n} \sum_{k=1}^{[nr]}g(c_{n} x_{k,n})\rightarrow
_{D}\omega
\ell_{G}(1,0),
\end{equation}
where $\ell_{G}(t,s)$ is the local time of the process $G(t)$ at
the spatial point $s$. We refer to \citet{BorIbr94} (and
their references for related work) for the particular situation where
$c_{n}x_{k,n} $ is a partial sum of i.i.d. random variables, and to
\citet{Ako93}, \citet{PhiPar98}, \citet{Jeg04} and
\citet{deJWan05} for the case where $c_{n}x_{k,n}$ is a partial
sum of a linear process. Wang and Phillips [(\citeyear{WanPhi09N1}), Theorem
2.1]
generalized these results to include not only linear process partial
sums but also cases where $c_{n}x_{k,n}$ is a partial sum of a
Gaussian process, including fractionally integrated time series.

Our present research on the statistic $S_{[nr]}$ in (\ref{tt5}) has a
similar motivation to this earlier work on convergence to a local time
process. However, the statistic $S_{[nr]}$ has a much more complex
U-statistic form, and the technical difficulties of establishing weak
convergence are greater. The approach of Wang and Phillips
[(\citeyear{WanPhi09N1}), Theorem 2.1] remains useful, however, and is
implemented in the proofs of Theorems \ref{th3}--\ref{th5}.

Finally we mention some earlier work investigating the intersection
local time process and weak convergence for certain specialized
situations. This work restricts the function $g$ in (\ref{tt5}) to the
indicator function and the discrete process $y_{k,n}$ in (\ref{mod1})
to a lattice random walk taking integer values; see, for instance,
\citet{Ald86}, \citet{vandenKon97}, \citet{vanKon01}
and \citet{vandenKon03}. The present
paper seems to the first to consider weak convergence to intersection
local time for a general linear process and a general
function $%
g$.

The proofs of Theorems \ref{th1a}--\ref{th3a} are given
in the supplementary material [\citet{WanPhi}].

\section{\texorpdfstring{Proofs of Theorems \protect\ref{th3}--\protect\ref{th5}}{Proofs of Theorems 3.1--3.3}}\label{sec6}

We start with several propositions. Their proofs are given in the
supplementary material [\citet{WanPhi}]. Throughout the section,
we let $C, C_1$, $C_2,\ldots$ be constants which may differ at each
appearance.
%
%
\begin{prop}\label{prop3}
$\!\!\!$Suppose Assumptions \ref{Assumption1} and \ref
{Assumption2} hold.
For any \mbox{$\al_1, \al_2\ge0$},
if ${\sup_x}|p(x)|<\infty$,
$\int(1+|x|^{\max\{[\al_1],[\al_2]\}+ 1})|p(x)|\,dx<\infty$
and $E|\ep_0|^{[\al_1]+[\al_2]+2}<\infty$, then
%
%
\begin{eqnarray} \label{st146}\qquad
\Lambda_n:\!&=& \mathop{\sum_{s,t=1}}_{s\not=t}^n
g(u_{s+1})g_1(u_{t+1})(1+|x_s|^{\al_1})(1+|x_t|^{\al_2})
p[(x_t-x_s)/h] \nonumber\\[-8pt]\\[-8pt]
&=& O_P (n^{3/2+\al_1/2+\al_2/2}h),\nonumber
\end{eqnarray}
where $g(x)$ and $g_1(x)$ are real functions such that
\[
\sup_{s\ge1} E\{[g^2(u_{s+1})+g_1^2(u_{s+1})] \mid{\cal
F}_s\}<\infty.
\]
If additionally $\al_1>0$, then
%
%
\begin{eqnarray} \label{st147}
\widetilde\Lambda_n:\!&=& \sum_{1\le s<t\le n}
g(u_{s+1})(1+|x_s|^{\al_1-1})p[(x_t-x_s)/h] \nonumber\\[-8pt]\\[-8pt]
&=& O_P \bigl(n^{\max\{3/2, 1+\al_1/2\}}h\bigr).\nonumber
\end{eqnarray}
\end{prop}
%
%
\begin{prop}\label{prop4} Suppose Assumptions \ref{Assumption1}--\ref
{Assumption3} hold. Then, for any
$g(x, \theta)$ satisfying (\ref{tz2}) and $|g(x, \theta)|\le
C(1+|x|^\beta)$, where $\theta\in\Omega_0$, we have
%
%
\begin{equation} \label{ho7}\qquad
\Delta_n :=\mathop{\sum_{s,t=1}}_{s\not= t}^nu_{s+1}g(x_t,
\theta)K[(x_{t}-x_s)/h] = O_P(n^{5/4+\beta/2}
h^{3/4}),
\end{equation}
provided that $nh^2\to\infty$, $nh^4\to0$, $\int(1+|x|^{\beta
+1})K(x)\,dx<\infty$ and $E|\ep_0|^{\beta+2}<\infty$.
Similarly, (\ref{ho7}) holds true if we replace
$g(x,\theta)$ and $\beta$ by $m(x)$ and $\nu$, respectively, where
$m(x)$ is defined as in Assumption \ref{Assumption6}.
\end{prop}
%
%
\begin{prop}\label{prop5} Suppose Assumptions \ref{Assumption1}--\ref
{Assumption3} hold and $nh^2\to
\infty$. Then, for any
real function $g(x)$ satisfying $\sup_{s\ge1} E\{g^2(u_{s+1})
\mid{\cal F}_s\}<\infty, $ we have
%
%
\begin{equation} \label{ho75}\qquad
\Gamma_n := \mathop{\sum_{s, t=1}}_{s\not=t}^n g(u_{s+1})(u_{t+1}^2-\si^2)
K^2[(x_{t}-x_s)/h] =o_P( n^{3/2}
h).
\end{equation}
\end{prop}

%
\begin{prop}\label{prop6} In addition to Assumptions \ref
{Assumption1}--\ref{Assumption3}, we have
$|u_j|\le
A$ and $nh^2\to\infty$. Then,
%
%
\begin{equation}\label{ho76}
R_n := \sum_{t=1}^{n}\mathop{\sum_{i, j=1}}_{
i\not= j}^{t-1}u_{i+1}u_{j+1}
K[(x_t-x_i)/h] K[(x_t-x_j)/h] = o_P( n^{3/2}
h).\hspace*{-30pt}
\end{equation}
\end{prop}
%
%
\begin{prop} \label{56}
Under Assumptions \ref{Assumption1}--\ref{Assumption3} and $h\log^2 n\to
0$, we have
%
%
\begin{equation} \label{fi2}
EZ_{tkr}^2 \le C \max_{1\le i,j\le n}E[|u_i|(1+|u_j|)]
\bigl(1+h\sqrt{t-r-k}\bigr)
\end{equation}
for $1\le k\le t-r$ and $r\ge1$,
where $Z_{tkr}= \sum_{i=k}^{t-r}u_{i+1}K[(x_t-x_i)/h]$.
Similarly,
%
%
\begin{equation}\label{d1}\qquad
E\Biggl\{\sum_{i=1}^{t-1}[u_{i+1}^2-E(u_{i+1}^2\mid{\cal
F}_j)]K^2[(x_t-x_i)/h] \Biggr\}^2\le
C \bigl(1+h\sqrt{t}\bigr).
\end{equation}
If in
addition $|u_j|\le A$, where $A$ is a constant, then
%
%
\begin{equation} \label{fi3}
EZ_{t12}^4 \le Ch^{3} t^{3/2},
\end{equation}
and for any
$1\le m\le t/2$,
%
%
\begin{equation} \label{m1}
EZ_{tm}^{*2} \le\frac
{Ch^2t^2}{m^{3/2}}+\frac{Ch^2t\log(t-m)}{ \sqrt m} +\frac
{Ch^2t}{m},
\end{equation}
where $Z_{tm}^{*}=\sum_{i=1}^{t-m-1}u_{i+1}
E(K[(x_t-x_i)/h]\mid{\cal F}_{t-m})$.
\end{prop}

%
\subsection{\texorpdfstring{Proof of Theorem \protect\ref{th3}}{Proof of Theorem 3.1}}\label{sec61}
By virtue of (\ref{s1}) and Theorem \ref{th5}, it suffices to verify
(\ref{lp1}) and show that
%
%
\begin{equation} \label{hp1}
S_{2n} =o_P\bigl(n^{3/4}\sqrt h\bigr) \quad\mbox{and}\quad
S_{3n}=o_P\bigl(n^{3/4}\sqrt h\bigr).
\end{equation}
To this end, for $
\delta>0$, let $\Omega_n=\{\hat{\theta}\dvtx
\|\hat{\theta}-\theta\|\le\delta\delta_n, \theta\in
\Omega_0\}$, where $\delta_n$ is given in Assumption \ref{Assumption4}(i).

We first prove (\ref{hp1}). Note that $\Omega_n\subset\Theta_0$ for
all $n$
sufficiently large. Under Assumption \ref{Assumption4},
it follows by Taylor's expansion that, whenever
$n$ is sufficiently large and $\hat\theta\in\Omega_n$,
%
%
\begin{equation} \label{t2}\qquad
S_{2n} = (\theta-\hat\theta) \mathop{\sum_{i, t=1}}_{i\not=t}^n
u_{i+1}\,
\frac{\partial f(x_t, \theta)}{\partial
\theta}K[(x_{t}-x_i)/h]+S_{2n1},
\end{equation}
where
\[
S_{2n1} \le C |\hat\theta-\theta|^2\mathop{\sum_{i, t=1}}_{i\not
=t}^n |u_{i+1}|
(1+ |x_t|^{\beta}) K[(x_{t}-x_i)/h].
\]
By Proposition \ref{prop4} with $g(x, \theta)=
\frac{\partial f(x, \theta)}{\partial\theta}$
and $\delta_n^2n^{1+\beta}\sqrt h\to0$,
the first term in the decomposition of $S_{1n}$ is equal to
\[
O_P(\delta_n n^{5/4+\beta/2}
h^{3/4}) = o_P\bigl(n^{3/4}\sqrt h\bigr).
\]
On the other hand, by Proposition \ref{prop3} and $nh^2\to\infty$,
we get
\[
S_{2n1}
=O_P(\delta_n^2 n^{3/2+\beta/2}h) = o_P\bigl(n^{3/4}\sqrt h\bigr).
\]
These facts imply, for any
$\delta>0$,
%
%
\begin{eqnarray} \label{mn1}
&&
P\bigl(|S_{2n}|\ge\delta n^{3/4}\sqrt h\bigr)\nonumber\\
&&\qquad\le P\bigl(|S_{2n}|\ge\delta n^{3/4}\sqrt h, \hat\theta\in
\Omega_n\bigr)+
P(\|\hat\theta-\theta\|\ge\delta\delta_n) \\
&&\qquad\to 0\qquad\mbox{as $n\to\infty$}.\nonumber
\end{eqnarray}
Similarly, by using Proposition \ref{prop3} and
noting
%
%
\begin{eqnarray} \label{hp2}
|S_{3n}| &\le&
C |\hat\theta-\theta|^2 \mathop{\sum_{i, t=1}}_{i\not=t}^n
\biggl|\frac{\partial f(x_i, \theta)}{\partial\theta}\biggr|
\biggl|\frac{\partial f(x_t, \theta)}{\partial\theta}\biggr|
K[(x_{t}-x_i)/h] \nonumber\\
&\le& C \delta_n^2 \mathop{\sum_{i, t=1}}_{i\not=t}^n
(1+|x_i|^{\beta}) (1+|x_t|^{\beta}) K[(x_{t}-x_i)/h]\\
&=& O_P(\delta_n^2
n^{3/2+\beta} h)=o_P\bigl(n^{3/4}\sqrt h\bigr),\nonumber
\end{eqnarray}
whenever
$\hat\theta\in\Omega_n$, we obtain,
for any $\delta>0$,
%
%
\begin{eqnarray} \label{mn2}
&&P\bigl(|S_{3n}|\ge\delta n^{3/4}\sqrt h\bigr)\nonumber\\
&&\qquad\le P\bigl(|S_{3n}|\ge\delta n^{3/4}\sqrt h, \hat\theta\in
\Omega_n\bigr)+
P(|\hat\theta-\theta|\ge\delta\delta_n) \\
&&\qquad\to 0\qquad\mbox{as $n\to\infty$}.\nonumber
\end{eqnarray}
Combining (\ref{mn1}) and (\ref{mn2}), we obtain (\ref{hp1}).

We next prove (\ref{lp1}). We may write
%
%
\begin{eqnarray} \label{mn3}
V_n^2 &=& \mathop{\sum_{s,t=1}}_{s\not=t}^n u_{s+1}^2u_{t+1}^2
K^2[(x_{t}-x_s)/h] \nonumber\\[-1pt]
&&{} +\mathop{\sum_{s,t=1}}_{s\not=t }^n
(\hat u_{s+1}^2-u_{s+1}^2)
\hat u_{t+1}^2 K^2[(x_{t}-x_s)/h]\nonumber\\[-9pt]\\[-9pt]
&&{} +
\mathop{\sum_{s,t=1}}_{s\not=t}^n u_{s+1}^2
(\hat u_{t+1}^2-u_{t+1}^2) K^2[(x_{t}-x_s)/h]\nonumber\\[-1pt]
:\!&=& V_{1n}+V_{2n}+V_{3n}.\nonumber
\end{eqnarray}
Recall
$|f(x_s, \theta)-f(x_s, \hat\theta)|\le
C \delta_n(1+|x_s|^{\beta})$ whenever $\hat\theta\in\Omega_n$ and
$|\hat u_{t+1}^2-u_{t+1}^2| =2|u_{t+1}||f(x_s, \theta)-f(x_s, \hat
\theta)|
+|f(x_s, \theta)-f(x_s, \hat\theta)|^2$. It is readily seen
from Proposition \ref{prop3} that, given $\hat\theta\in\Omega_n$,
\begin{eqnarray*}
|V_{2n}|+ |V_{3n}|&\le& C\delta_n \mathop{\sum_{s,t=1}}_{s\not=t}^n
|u_{s+1}| u_{t+1}^2 (1+|x_s|^{\beta})
K[(x_{t}-x_s)/h] \\
&&{}+C\delta_n^2 \mathop{\sum_{s,t=1}}_{s\not=t}^n
u_{t+1}^2 (1+|x_s|^{2\beta})
K[(x_{t}-x_s)/h]\\
&&{}+C\delta_n^3 \mathop{\sum_{s,t=1}}_{s\not=t}^n
|u_{s+1}| (1+|x_s|^{\beta}) (1+|x_t|^{2\beta})
K[(x_{t}-x_s)/h]\\
&&{}+C\delta_n^4 \mathop{\sum_{s,t=1}}_{s\not=t }^n
(1+|x_s|^{2\beta}) (1+|x_t|^{2\beta})
K[(x_{t}-x_s)/h]\\
&=&O_P(n^{3/2}h)
(\delta_n n^{\beta/2}+
\delta_n^2 n^{\beta}+\delta_n^3 n^{3\beta/2}+\delta_n^4 n^{2\beta
})\\
&=& o_P(n^{3/2}h),\vadjust{\goodbreak}
\end{eqnarray*}
since $nh^2\to\infty$ and $\delta_n^2n^{1+\beta}\sqrt h\to0$.
As for $V_{1n}$, by Proposition \ref{prop5},
we have
\begin{eqnarray*}
V_{1n}
&=& \si^4 \mathop{\sum_{s,t=1}}_{s\not=t}^n
K^{2}[(x_{t}-x_s)/h] +
\mathop{\sum_{s,t=1}}_{s\not=t}^n (u_{t+1}^2+\si^2)
(u_{s+1}^2-\si^2)
K^{2}[(x_{t}-x_s)/h] \\
&=& \si^4 \mathop{\sum_{s,t=1}}_{s\not=t}^n
K^{2}[(x_{t}-x_s)/h] +o_P(n^{3/2}h).
\end{eqnarray*}
Taking these estimates into (\ref{mn3}), we get
the first part of (\ref{lp1}).

In order to prove the second part of (\ref{lp1}), we first assume
$|u_j|\le A$.
In this case, simple calculations together with Propositions \ref{prop5}
and \ref{prop6} yield that
%
%
\begin{eqnarray} \label{88}\qquad
\sum_{t=2}^{n}Y_{nt}^{2} &=& \sum_{t=2}^{n}
\sum_{s=1}^{t-1}u_{s+1}^2K^2[(x_{t}-x_s)/h] \nonumber\\
&&{}+\sum
_{t=1}^{n}\mathop{\sum_{i, j=1}}_{i\not= j}^{t-1}u_{i+1}u_{j+1}
K[(x_t-x_i)/h] K[(x_t-x_j)/h] \\
&=&\frac{\si^2}2\mathop{\sum_{s,t=1}}_{s\not=t}^{n} K^2
[(x_{t}-x_s)/h]
+o_P(n^{3/2}h)\nonumber
\end{eqnarray}
as required. The idea to remove the restriction $|u_j|\le A$ is the
same as
in the proof of Theorem \ref{th5}. We omit the details.
The proof of Theorem \ref{th3} is now complete.

%
\subsection{\texorpdfstring{Proof of Theorem \protect\ref{th4}}{Proof of Theorem 3.2}}\label{sec62}

Put $\hat{u}_{t+1}^*=u_{t+1}+f(x_t, \theta)-f(x_{t},
\hat{\theta})$. Under~$H_1$, we may write
%
%
\begin{equation} \label{202}
S_n = S_{1n}+2S_{2n}+S_{3n} -S_{4n} +S_{5n},
\end{equation}
where $S_{1n}, S_{2n}, S_{3n}$ are defined as in (\ref{s1}), and
\begin{eqnarray*}
S_{4n}&=&2\rho_n \mathop{\sum_{i, t=1}}_{i\not=t}^n m(x_i) \hat
{u}_{t+1}^* K[(x_{t}-x_i)/h], \\
S_{5n} &=&\rho_n^2 \mathop{\sum_{i, t=1}}_{i\not=t}^n m(x_i) m(x_t)
K[(x_{t}-x_i)/h].
\end{eqnarray*}
Thus (\ref{205}) will follow if we prove
%
%
\begin{eqnarray}
\label{206}
S_{jn} & =& O_P(n^{3/4}h^{1/2}),\qquad j=1, 2, 3, \\
\label{206a}
S_{4n}& =& O_P(\rho_nn^{5/4+\nu/2}h^{3/4}), \\
\label{209}
V_n^2 &=& O_P(n^{3/2}h+\rho_n^4 n^{3/2+2\nu}h)\qquad \mbox{under
$H_1$,}
\end{eqnarray}
and for any $\ep_n\to0$,
%
%
\begin{equation}\label{208}
S_{5n} \ge\ep_n \rho_n^2n^{3/2+\nu}h\qquad \mbox{in Probab.}
\end{equation}
Here and below, the notation $A_n\ge B_n$, in Probab. means that
$\lim_{n\to\infty}\!P(A_n\ge B_n)=1$, as $n\to\infty$.
Indeed, by choosing $\ep_n^{-2}=\min\{\rho_n^2n^{1/2+\nu}\sqrt h,
n^{3/2}\sqrt h\}$,
it is readily seen that $\ep_n\to0$,
$|S_{jn}|=O_P(\ep_n S_{5n})=o_P(S_{5n})$ for $j=1,2,3,4$ and
$S_{5n}/V_n\ge\ep_n^{-1}$, in Probab.
Hence $S_n/V_n\ge\ep_n^{-1}/2$, in Probab., which yields~(\ref{205}).

We next prove (\ref{206a})--(\ref{208}). The proof of (\ref{206})
for $j=2,3$ is given in (\ref{hp1}), and
the result for $j=1$ is simple by martingale properties and Proposition
\ref{56}.

Equation (\ref{208}) first. We may write
%
%
\begin{equation}
S_{5n} = S_{5n1}+ S_{5n2},
\end{equation}
where $S_{5n1}=2\rho_n^2 \sum_{1\le i<t\le n} m^2(x_i)
K[(x_{t}-x_i)/h]$ and
\[
|S_{5n2}| \le2\rho_n^2 \sum_{1\le i<t\le n} |m(x_i)|
|m(x_t)-m(x_i)|
K[(x_{t}-x_i)/h].
\]
Let $\nu'=\nu$ if $\nu>0$ and $\nu'=\gamma'$ if $\nu=0$. It
follows from
(\ref{2011}) and Proposition~\ref{prop3} that
%
%
\begin{eqnarray} \label{us1}
|S_{5n2}|
&\le&C h^{\gamma} \rho_n^2 \sum_{1\le i<t\le n} (1+|x_i|^{\nu})
(1+|x_i|^{\nu'-1}+|x_t-x_i|^{\nu})
K_{\gamma}[(x_{t}-x_i)/h] \nonumber\\
&\le& C h^{\gamma} \rho_n^2 \sum_{1\le i<t\le n}
\{(1+ |x_i|^{\nu'-1}+ |x_i|^{\nu+\nu'-1}) K_\gamma
[(x_{t}-x_i)/h] \nonumber\\
&&\hspace*{107pt}{} +h^{\nu}
(1+|x_i|^{\nu}) K_{\nu+\gamma}[(x_{t}-x_i)/h]\} \\
&= & O_P\bigl( h^{1+\gamma} \rho_n^2 \bigl[n^{\max\{3/2, 1+(\nu+\nu')/2\}
}+n^{3/2+\nu/2}\bigr]\bigr)\nonumber\\
&=&
O_P( h^{1+\gamma} \rho_n^2 n^{3/2+\nu}),\nonumber
\end{eqnarray}
where $K_u(x)=|x|^uK(x), u>0$ and we have used the fact that ${\sup
_x}|K_u(x)|<\infty$ whenever
$\int K_u(x) \,dx<\infty$ [recall $\sup_x|K(x)|<\infty$].
Since $h\to0$ and $0<\gamma\le1$, to prove (\ref{208}), it only
needs to show that,
for any $h^{\gamma/2}\le\ep_n\to0$,
%
%
\begin{equation} \label{208a}
S_{5n1}\ge\ep_n \rho_n^2 n^{3/2+\nu}\qquad \mbox{in
Probab.}
\end{equation}
In fact,
by (\ref{no11}) and letting $x_{[ns], n}=x_{[ns]}/(\sqrt n\phi)$,
%
%
\begin{eqnarray} \label{er1}\qquad
\inf_{n/2\le j\le n}|x_j|&\ge& \sqrt n\phi
\Bigl( \inf_{1/2\le s\le1}|G_1(s)|-\sup_{1/2\le s\le
1}\bigl|x_{[ns],n}-G_1(s)\bigr|\Bigr)\nonumber\\[-8pt]\\[-8pt]
&\ge& \ep_n^{1/4\nu} \sqrt n\qquad \mbox{in Probab.}\nonumber
\end{eqnarray}
Similarly, by using
(\ref{r28}) in Theorem \ref{th1a}, we have
%
%
\begin{equation} \label{er2}
\sum_{n\ge t>i\ge n/2} K[(x_{t}-x_i)/h] \ge\ep_n^{1/4}
n^{3/2} h\qquad
\mbox{in Probab.}
\end{equation}
Combining (\ref{201}), (\ref{er1}) and (\ref{er2}), we obtain that
\begin{eqnarray*}
S_{5n1} &\ge& \ep_n^{1/2} \rho_n^2 \sum_{n\ge t>i\ge n/2}
|x_i|^{2\nu}I(|x_i|\ge1) K[(x_{t}-x_i)/h] \\
&\ge& \ep_n^{3/4} \rho_n^2 n^{\nu} \sum_{n\ge t>i\ge n/2}
K[(x_{t}-x_i)/h] \\
&\ge& \ep_n \rho_n^2 n^{3/2+\nu} h\qquad \mbox{in Probab.}
\end{eqnarray*}
This provides (\ref{208a}) and also completes the proof of (\ref{208}).

Next prove (\ref{206a}). We have
\[
S_{4n}=
2\rho_n \mathop{\sum_{i, t=1}}_{i\not=t}^n m(x_i) u_{t+1} K
[(x_{t}-x_i)/h] +S_{4n1},
\]
where, by recalling $|f(x_t, \theta)-f(x_{t},\hat{\theta})|\le C
\|\hat\theta-\theta\| (1+|x_t|^{\beta})$ by Assumption~\ref
{Assumption4}, it
follows from Proposition \ref{prop3} that
\begin{eqnarray*}
|S_{4n1}| &\le&
\mathop{\sum_{i, t=1}}_{i\not=t}^n(1+ |x_i|^{\nu})|f(x_t, \theta)-f(x_{t},
\hat{\theta})| K[(x_{t}-x_i)/h] \\
&\le& C\rho_n \|\hat\theta-\theta\|
\mathop{\sum_{i, t=1}}_{i\not=t}^n (1+|x_i|^{\nu}) (1+|x_t|^{\beta})
K[(x_{t}-x_i)/h] \\
&=& O_P(\rho_n \delta_n n^{3/2+\nu/2+\beta/2}h).
\end{eqnarray*}
This, together with Proposition \ref{prop4}, yields that
\begin{eqnarray*}
S_{4n} &=& O_P(\rho_n \delta_n n^{3/2+\nu/2+\beta/2}h)+O_P
( \rho_n n^{5/4+\nu/2}
h^{3/4}) \\
&=&O_P( \rho_n n^{5/4+\nu/2}
h^{3/4}),
\end{eqnarray*}
since $\delta_n^2n^{1+\beta}\sqrt n\to0$. The result (\ref{206a})
is proved.

Finally, we prove (\ref{209}). Under $H_1$, we have
%
%
\begin{eqnarray} \label{203}\qquad
V_n^2 &=& \sum_{s,t=1,s\not=t}^{n}[\hat{u}_{t+1}^*+\rho_n
m(x_t)]^{2}
[\hat{u}_{s+1}^*+\rho_n m(x_s)]^{2} K^{2}[(x_{t}-x_{s})/h]
\nonumber\\[-8pt]\\[-8pt]
&\le& 2 V_{6n} +4 V_{7n} +2 V_{8n},\nonumber
\end{eqnarray}
where
\begin{eqnarray*}
V_{6n}&=& \mathop{\sum_{s, t=1}}_{s\not=t}^n\hat{u}_{t+1}^{*2}
\hat{u}_{s+1}^{*2} K^2[(x_{t}-x_s)/h],\\
V_{7n} &=&\rho_n^2 \mathop{\sum_{s, t=1}}_{s\not=t}^n \hat{u}_{t+1}^{*2}m^2(x_s)
K^2[(x_{t}-x_s)/h], \\
V_{8n} &=&\rho_n^4 \mathop{\sum_{s, t=1}}_{s\not=t}^n m^2(x_t)
m^2(x_s)K^2[(x_{t}-x_s)/h].
\end{eqnarray*}
By recalling $|m(x)|\le C |x|^{\nu}$ and
\begin{eqnarray*}
\hat{u}_{t+1}^{*2}&\le&2\bigl(u_{t+1}^2+
|f(x_t, \theta)-f(x_{t},\hat{\theta})|^2\bigr)\\
&\le& C
[u_{t+1}^2+O_P(\delta_n^2)(1+|x_t|^{2\beta})],
\end{eqnarray*}
it following repeatedly from Proposition \ref{prop3} and $\delta
_n^2n^{1+\beta}\sqrt h\to0$ that
\begin{eqnarray*}
V_{6n}
&\le& C\mathop{\sum_{s, t=1}}_{s\not=t}^n[u_{s+1}^2+
O_P(\delta_n^2)(1+|x_s|^{2\beta})] [u_{t+1}^2+O_P(\delta_n^2)
(1+|x_t|^{2\beta})]\\
&&\hspace*{28pt}{}\times K^2[(x_{t}-x_s)/h]\\
&=&O_P(n^{3/2}h)+O_P(\delta_n^2n^{3/2+\beta}h)+O_P(\delta
_n^4n^{3/2+2\beta}h)\\
&=& O_P(n^{3/2}h).
\end{eqnarray*}
Similarly, we have
\begin{eqnarray*}
V_{7n} &\le& C\rho_n^2 \mathop{\sum_{s, t=1}}_{s\not=t}^n[u_{s+1}^2+
O_P(\delta_n^2)(1+|x_s|^{2\beta})]
(1+|x_t|^{2\nu})] K^2[(x_{t}-x_s)/h]\\
&=&O_P(\rho_n^2 n^{3/2+\nu}h)+O_P(\rho_n^2 \delta_n^2
n^{3/2+\beta+\nu}h) =O_P(\rho_n^2 n^{3/2+\nu}h),\\
V_{8n} &\le&\rho_n^4 \mathop{\sum_{s, t=1}}_{s\not=t}^n (1+|x_t|^{2\nu
})
(1+|x_s|^{2\nu})K^2[(x_{t}-x_s)/h]\\
&=&O_P (\rho_n^4 n^{3/2+2\nu}h).
\end{eqnarray*}
Combining all these estimates, we obtain
\begin{eqnarray*}
V_n^2&=&O_P(n^{3/2}h)+O_P(\rho_n^2 n^{3/2+\nu}h)+O_P (\rho_n^4
n^{3/2+2\nu}h) \\
&=&O_P(n^{3/2}h+\rho_n^4 n^{3/2+2\nu}h)
\end{eqnarray*}
as required. The proof of Theorem \ref{th4} is complete.

%
\subsection{\texorpdfstring{Proof of Theorem \protect\ref{th5}}{Proof of Theorem 3.3}}\label{sec6.3}

We first assume $|u_t|\le A$, where $A$ is a constant. This restriction
will be removed later. Write $G_n(t)=x_{[nt]}/\sqrt n\phi$ and
$V_n(t)=\sum_{j=1}^{[nt]}u_{j+1}/\sqrt n \si$. Under Assumptions
\ref{Assumption1} and \ref{Assumption2}, the same arguments as those in
\citet{BucCha07} or \citet{WanPhi09N2}, with minor
modifications, show that
%
%
\begin{equation} \label{ne1}
(G_{n},V_{n})\Rightarrow_{D}(G,V)
\end{equation}
on $D[0,1]^2$,
where $G(t)= W(t)+\kappa\int_0^te^{\kappa(t-s)} W(s)\,ds$. By
virtue of (\ref{ne1}),
it follows from the so-called
Skorohod--Dudley--Wichura representation theorem that there is a
common
probability space $(\Omega,\mathcal{F},P)$ supporting $%
(G_{n}^{0},V_{n}^{0}) $ and $(G,V)$ such that
%
%
\begin{equation} \label{p2}
(G_{n},V_{n})=_{d}(G_{n}^{0},V_{n}^{0}) \quad\mbox{and}\quad
(G_{n}^{0},V_{n}^{0})\rightarrow_{\mathrm{a.s.}}(G,V)
\end{equation}
in $D[0,1]^{2}$ with the uniform topology. Moreover, as in the proof
of Lemma~2.1 in \citet{ParPhi01}, $V_{n}^{0}$ can be chosen
such that, for each $n\geq1$,
%
%
\begin{equation} \label{p3}
V_{n}^{0}(k/n)=V(\tau_{nk}/n),\qquad k=1,2,\ldots,n,
\end{equation}
where $\tau_{n, k},1\le k\le n$, are stopping times with respect to
${\cal F}_{n, k}^0$ in $(\Omega,%
\mathcal{F},P)$ with
\[
{\cal F}_{n, k}^0 =\si\{ V(r), r\le\tau_{n, k}/n;
G_{n}^{0}(s/n), s=1,\ldots,k+1\},
\]
satisfying $\tau_{n,0}=0$,
%
%
\begin{equation} \label{p4}
\sup_{1\leq k\leq n}\biggl|\frac{\tau_{n, k}-k}{n^{\delta}}\biggr|
\rightarrow_{\mathrm{a.s.}}0
\end{equation}
as $n\rightarrow\infty$ for any $1/2<\delta<1$, and
%
%
\begin{eqnarray}\label{p4a}
E [(\tau_{n, k}-\tau_{n,k-1})\mid{\cal F}_{n, k-1}^0]
&=&\si^{-2} E [u_{k+1}^2\mid{\cal F}_{k}] \quad\mbox
{and} \nonumber\\[-8pt]\\[-8pt]
E [(\tau_{n, k}-\tau_{n,k-1})^{2m}\mid{\cal F}_{n,
k-1}^0] &\le& C \si^{-4m} E [u_{k+1}^{4m}\mid{\cal
F}_{k}],\quad m\ge1\mbox{,\quad a.s.}\hspace*{-40pt}
\nonumber
\end{eqnarray}
for some constant $C>0$. We mention that result (\ref{p4a}) does
not explicitly appear in Lemma 2.1 of \citet{ParPhi01};
however, it can be obtained by a construction along the same lines as
Theorem A1 of \citet{HalHey80}.\vadjust{\goodbreak}

It follows from (\ref{p3}) that, under the extended probability
space,
%
%
\begin{eqnarray} \label{mo100}
&& \Biggl(\frac1{\si d_n}\sum_{t=2}^nu_{t+1} Y_{nt},
\frac1{d_n^2}\sum_{t=2}^n Y_{nt}^2\Biggr) \nonumber\\[-8pt]\\[-8pt]
& &\qquad =_{d} \Biggl( \sum_{t=2}^n [V(\tau_{n,t}/n)-V(\tau
_{n,t-1}/n)] Y_{n, t}^* , \frac1{n} \sum_{t=2}^n
Y_{nt}^{*2}\Biggr),\nonumber
\end{eqnarray}
where, with $c_n=\sqrt n\phi/h$,
\[
Y_{nt}^*=\frac{ n\si}{d_n}
\sum_{i=1}^{t-1}[V(\tau_{n,i}/n)-V(\tau_{n,i-1}/n)]
K\{c_n [G_{n}^{0}(t/n)-G_{n}^{0}(i/n)]\}.
\]

To establish our main result, we extend $\sum_{i=2}^n [V(\tau
_{n,t}/n)-V(\tau_{n,t-1}/n)] Y_{n, t}^*$ to a continuous
martingale. This can be done by defining
%
%
\begin{equation} \label{tt30}
M_{n}(r)\,{=}\,\sum_{t=2}^{j-1} Y_{nt}^*\biggl[\!V\!\biggl(\frac{\tau
_{n,t}}{n}\!\biggr)\,{-}\,V\!\biggl(\frac{\tau
_{n,t-1}}{n}\biggr)\!\biggr]\,{+}\,
Y_{n,j}^*\biggl[\!V(r)\,{-}\,V\!\biggl(\frac{\tau_{n,j-1}}{n}\biggr)\!\biggr]\hspace*{-45pt}
\end{equation}
for $\tau_{n,j-1}/n<r\leq\tau_{n,j}/n,j=1,2,\ldots,n$, and
%
%
\begin{equation} \label{tt30a}
M_{n}(r)\,{=}\,\sum_{t=2}^{n} Y_{nt}^*\biggl[\!V\!\biggl(\frac{\tau
_{n,t}}{n}\biggr)\,{-}\,V\!\biggl(\frac{\tau
_{n,t-1}}{n}\biggr)\!\biggr]\,{+}\,\frac1{\sqrt n}\biggl[\!V(r)\,{-}\,V\!\biggl(\frac{\tau
_{n,n}}{n}\biggr)\!\biggr]\hspace*{-45pt}
\end{equation}
for $r\ge\tau_{n,n}/n$. It is readily seen that $M_{n}$ is a
continuous martingale with quadratic variation process $[M_{n}]$
given by
%
%
\begin{equation} \label{pw1}
[ M_{n}]_{r} = \sum_{t=2}^{j-1}Y_{nt}^{*2} \biggl(\frac{\tau
_{n,t}}{n}-%
\frac{\tau_{n,t-1}}{n}\biggr)+ Y_{n,j}^{*2} \biggl(r-\frac{\tau
_{n,j-1}}{n}\biggr)
\end{equation}
for $\tau_{n,j-1}/n<r\leq\tau_{n,j}/n,j=1,2,\ldots,n$, and
%
%
\begin{equation} \label{pw11}
[ M_{n}]_{r} = \sum_{t=2}^{n}Y_{nt}^{*2} \biggl(\frac{\tau
_{n,t}}{n}-%
\frac{\tau_{n,t-1}}{n}\biggr)+ \frac1{ n} \biggl(r-\frac{\tau
_{n,n}}{n}\biggr)
\end{equation}
for $r\ge\tau_{n,n}/n$. Similarly, the covariance process
$[M_{n},V]$ of $M_{n}$ and $V$ is given by
%
%
\begin{equation} \label{pwa1}
[ M_{n},V]_{r} = \sum_{t=2}^{j-1}Y_{nt}^* \biggl(\frac{\tau
_{n,t}}{n}-%
\frac{\tau_{n,t-1}}{n}\biggr) + Y_{n,j}^* \biggl(r-\frac{\tau
_{n,j-1}}{n}\biggr)
\end{equation}
for $\tau_{n,j-1}/n<r\leq\tau_{n,j}/n,j=1,2,\ldots,n$, and
%
%
\begin{equation} \label{pwa1a}
[ M_{n},V]_{r} = \sum_{t=2}^{n}Y_{nt}^* \biggl(\frac{\tau
_{n,t}}{n}-%
\frac{\tau_{n,t-1}}{n}\biggr) + \frac1{\sqrt n} \biggl(r-\frac
{\tau
_{n,n}}{n}\biggr)
\end{equation}
for $r\ge\tau_{n,n}/n$.\vadjust{\goodbreak}

Write $\rho_{n}(t)=\inf\{s\dvtx[M_{n}]_{s}>t\}$, a
sequence of time changes. Note that $[M_n]_{\infty}=\infty$
for every $n\ge1$ and
%
%
\begin{equation} \label{tt31}
[ M_{n},V]_{\rho_{n}(t)}\rightarrow_{P}0\qquad \mbox{as
$n\to\infty$}
\end{equation}
for every $t\in R$, by (\ref
{fin2}) in Proposition \ref{prop2a} below.
Theorem
2.3 of Revuz and Yor [(\citeyear{RevYor99}), page 524]
yields that, if we call $B^{n}$ [i.e.,
$B^{n}(r)=M_{n}\{\rho_{n}(r)\}$] the DDS Brownian motion [see, e.g.,
\citet{RevYor99}, page 181] of the continuous
martingale $M_{n}$ defined by (\ref{tt30}) and (\ref{tt30a}),
then $B^{n}$ converges in
distribution to a Wiener process $W$. Since the law of the processes
$B^{n}$ are all given by Wiener measure,
it is plain that $B^n(r)\Rightarrow W(r)$ (mixing),
where the concept of mixing can be found in
Hall and Heyde (\citeyear{HalHey80}), page 56.
This, together with (\ref{fin4}) in Proposition \ref{prop2} below,
yields that $(B^n(r), [M_{n}]_{1}) \Rightarrow(W(r), \eta^2)$, where
$W$ is independent of
$\eta^2= L_{G}(1,0)$, defined as in (\ref{tt2}).
Now, by noting that $M_{n}(1)$ is equal to $B^{n}([M_{n}]_{1})$,
the continuous mapping theorem implies that
%
%
\begin{equation}\label{tt32}
(M_{n}(1),[M_{n}]_{1})\rightarrow_{D}(\eta N,\eta^{2}),
\end{equation}
where $N$ is a normal variate independent of $\eta$.

By virtue of (\ref{mo100}) and (\ref{tt32}), the required result
of the theorem follows~(\ref{pw19}) and~(\ref{pw31}) in
Proposition \ref{prop2} and Proposition \ref{prop3a} below.

It remains to show the following Propositions \ref{prop2a}--\ref
{prop3a}, whose proofs are given in the supplementary material
[\citet{WanPhi}]. The proof of Theorem~\ref{th5} under $|u_j|\le A
$ is now complete.
%
%
\begin{prop} \label{prop2a}
In addition to Assumptions \ref{Assumption1}--\ref{Assumption3}, assume
that \mbox{$|u_j|\le A$},
$nh^2\to\infty$ and $h\log^2n\to0$. Then, as $n\to\infty$,
%
%
\begin{equation} \label{fin2}
[M_n, V]_r \to0\qquad \mbox{in Probab.}
\end{equation}
uniformly
on $r\in[0, T]$, where $T$ is an arbitrary given constant.
\end{prop}
%
%
\begin{prop} \label{prop2}
In addition to Assumptions \ref{Assumption1}--\ref{Assumption3}, assume
that \mbox{$|u_j|\le A$}, $nh^2\to
\infty$ and $nh^4\log^2 n\to0$. Under the extended probability
space used in~(\ref{p2}),
we have
%
%
\begin{equation}\label{fin4}
[M_n]_1 \to_P \eta^2,
\end{equation}
where $\eta^2=
L_{G}(1,0)$ is defined as in (\ref{tt2}), and
%
%
\begin{equation}\label{pw19}
[M_n]_1-\frac1{n}
\sum_{t=1}^nY_{nt}^{*2} =o_P(1).
\end{equation}
\end{prop}
%
%
\begin{prop} \label{prop3a}
In addition to Assumptions \ref{Assumption1}--\ref{Assumption3}, assume
that \mbox{$|u_j|\le A$}, $nh^2\to
\infty$ and $nh^4\log^2 n\to0$. Then,
%
%
\begin{equation}\label{pw31}
M_{n}(1)-\sum_{t=2}^{n} Y_{nt}^*\biggl[V\biggl(\frac{\tau
_{n,t}}{n}\biggr)-V\biggl(\frac{\tau_{n,t-1}}{n}\biggr)\biggr] = o_P(1).
\end{equation}
\end{prop}

We next remove the restriction $|u_j|\le A$. To this end, let
\begin{eqnarray*}
u_{1j} &=& u_jI(|u_j|\le A/2)-E[u_jI(|u_j|\le A/2)\mid{\cal
F}_{j-1}], \\
u_{2j} &=& u_jI(|u_j|> A/2)-E[u_jI(|u_j|> A/2)\mid{\cal
F}_{j-1}]
\end{eqnarray*}
and
\[
Y_{1nt} = \sum_{i=1}^{t-1} u_{1,i+1} K[(x_t-x_i)/h],\qquad
Y_{2nt} = \sum_{i=1}^{t-1} u_{2,i+1} K[(x_t-x_i)/h].
\]
With this notation, we may write
%
%
\begin{eqnarray}\label{find1}
\frac1{d_n}\!\sum_{t=2}^nu_{t+1} Y_{nt} &=& \frac1{d_n}\!\sum_{t=2}^nu_{1,
t+1} Y_{1nt}\,{+}\,\frac1{d_n}\!\sum_{t=2}^nu_{1, t+1} Y_{2nt}\,{+}\,\frac1{d_n}\!\sum
_{t=2}^nu_{2, t+1} Y_{nt}\nonumber\hspace*{-40pt}\\[-8pt]\\[-8pt]
:\!&=& \frac1{d_n}\!\sum_{t=2}^nu_{1, t+1} Y_{1nt}\,{+}\,\Lam_{1n}\,{+}\,\Lam_{2n},
\nonumber\hspace*{-40pt}\\\label{find21}
\frac1{d_n^2}\!\sum_{t=2}^n Y_{nt}^2 &=& \frac1{d_n^2}\!\sum_{t=2}^n
Y_{1nt}^2\,{+}\,\frac2{d_n^2}\!\sum_{t=2}^n Y_{1nt} Y_{2nt}\,{+}\,\frac1{d_n^2}\!\sum
_{t=2}^n Y_{2nt}^2 \nonumber\hspace*{-40pt}\\[-8pt]\\[-8pt]
:\!&=& \frac1{d_n^2}\!\sum_{t=2}^n Y_{1nt}^2\,{+}\,\Lam_{3n}\,{+}\,\Lam_{4n}.
\nonumber\hspace*{-40pt}
\end{eqnarray}
Recall that $|u_{1j}|\le A$, and $u_{1j}$ is a martingale
difference satisfying
\begin{eqnarray*}
E(u_{1t}^2\mid{\cal F}_{t-1}) &=&
E\bigl(u_{t}^2I(|u_t|\le A)\mid{\cal F}_{t-1}\bigr)\\
&&{}- \bigl[E\bigl(u_{t}I(|u_t|\le
A)\mid{\cal F}_{t-1}\bigr)\bigr]^2\\
&\to&\si^2\qquad \mbox{a.s.}
\end{eqnarray*}
as $j,
A\to\infty$.
It follows from the proof of (\ref{192}) under $|u_j|\le A$ that, when
$n\to\infty$ first, and then $A\to\infty$,
%
%
\begin{equation}
\Biggl(\frac1{\si d_{n}}\sum_{t=2}^nu_{1,t+1} Y_{1nt},
\frac1{d_{n}^2}\sum_{t=2}^n Y_{1nt}^2\Biggr) \to_D (\eta N, \eta^2).
\end{equation}
Now it is readily seen that the required result will follow if
we prove
%
%
\begin{equation}\label{pro}
\Lambda_{in} \to_P 0,\qquad i=1,2,3,4,
\end{equation}
as
$n\to\infty$ first, and then $A\to\infty$. In fact, by virtue of
(\ref{fi2}) in Proposition~\ref{56},
\[
\sup_{1\le i\le n} Eu_i^2\le\sup_{1\le i\le n}(Eu_i^4)^{1/4}<\infty
\]
and $\sup_xK(x)<\infty$, we have, for $1\le t\le n$,
\begin{eqnarray*}
EY_{nt}^2 &\le& 2\sup_xK(x) Eu_{t}^2+2
E\Biggl(\sum_{i=1}^{t-2} u_{i+1} K[(x_t-x_i)/h]\Biggr)^2 \\
&\le& C \sup_{1\le i\le n} Eu_{i}^2 \bigl(1+h^2\sqrt t\log t+h\sqrt
t\bigr)\le C_1h\sqrt n,
\end{eqnarray*}
since $h\log n\to0$ and
$nh^2\to\infty$. Similarly,
\begin{eqnarray*}
EY_{1nt}^2 &\le& C \sup_{1\le
i\le n}Eu_{i}^2I(|u_i|\le A)
\bigl(1+h^2\sqrt t\log t+h\sqrt t\bigr) \le C_1h\sqrt n, \\
EY_{2nt}^2 &\le& C \sup_{1\le i\le n}Eu_{i}^2I(|u_i|> A)
\bigl(1+h^2\sqrt t\log t+h\sqrt t\bigr) \le C_1A^{-2} h\sqrt n.
\end{eqnarray*}
These results, together with the fact that $u_{1j}$ and
$u_{2j}$ both are martingale difference satisfying
\begin{eqnarray*}
\sup_j
E(u_{1,j+1}^2\mid{\cal F}_j) &\le&
\sup_j [E (u_j^4\mid{\cal F}_j)]^{1/2}\le C,\\
\sup_j E(u_{2,j+1}^2\mid{\cal F}_j) &\le&
\sup_j E \bigl(u_j^2I_{|u_j|>A}\mid{\cal F}_j\bigr)\\
&\le& A^{-2}\sup_{j}
E (u_j^4\mid{\cal F}_j) \le CA^{-2},
\end{eqnarray*}
yield that, as $n\to\infty$ first, and then $A\to\infty$,
\begin{eqnarray*}
E\Lam_{1n}^2 &\le& \frac C{n^{3/2}h}
\sum_{t=2}^n EY_{2nt}^2 \le CA^{-2}\to0,\\
E\Lam_{2n}^2 &\le& \frac{CA^{-2}}{n^{3/2}h}
\sum_{t=2}^n EY_{nt}^2 \le CA^{-2}\to0,\\
E\Lam_{4n} &\le& \frac{C}{n^{3/2}h}
\sum_{t=2}^n EY_{2nt}^2 \le CA^{-2}\to0,\\
E|\Lam_{3n}| &\le& \frac{C}{n^{3/2}h} \sum_{t=2}^n
(EY_{1nt}^2)^{1/2}(EY_{2nt}^2)^{1/2} \le CA^{-1}\to0.
\end{eqnarray*}
This
proves (\ref{pro}), and hence the proof of Theorem \ref{th5} is
complete.

\section*{Acknowledgments}

Our thanks to the Editor, Associate Editor and referees for helpful
comments on earlier versions.

\begin{supplement}
\stitle{Supplement to ``A specification test for nonlinear nonstationary models''}
\slink[doi]{10.1214/12-AOS975SUPP} 
\sdatatype{.pdf}
\sfilename{aos975\_supp.pdf}
\sdescription{Further details on the derivations in\vadjust{\goodbreak} the present paper
and supporting lemmas and proofs of the main results on convergence to
intersection local time are contained in the supplement to the paper,
\citet{WanPhi}.}
\end{supplement}

%

\printaddresses

\end{document}